\documentclass[italian,a4paper,10pt]{article}
\usepackage[latin1]{inputenc}
\usepackage[english]{babel}
\usepackage{babel,amsmath,amssymb,amsbsy,amsfonts,latexsym,amsthm,mathrsfs}
\usepackage{xcolor}
\usepackage{version}
\usepackage{verbatim}
\usepackage[title]{appendix}

\usepackage{amsmath,amsfonts,latexsym}
\newtheorem{theorem}{Theorem}[section]
\newtheorem{proposition}[theorem]{Proposition}
\newtheorem{definition}[theorem]{Definition}
\newtheorem{lemma}[theorem]{Lemma}
\newtheorem{corollary}[theorem]{Corollary}
\newtheorem{remark}[theorem]{Remark}
\newtheorem{example}[theorem]{Example}
\newtheorem{assumption}[theorem]{Assumption}

\def\proof{{\sl Proof. }}

\def\cl#1{{\mathscr #1}}

\def\P{{\mathbb P}}
\def\R{{\mathbb R}}

\def\N{{\mathbb  N}}

\def\E{{\mathbb E}}

\def\Var{{\rm Var }}
\def\Cov{{\rm Cov }}

\def\<{\langle}
\def\>{\rangle}

\newcommand{\cvd}{\hfill \ensuremath{\Box}\medskip}
\def\appendix{\par
  \setcounter{chapter}{0}
  \def\@chapter{Appendix}
  \def\thechapter{\Alph{chapter}}}
\font\tengoth=eufm10 scaled 1200 \font\sevengoth=eufm7 scaled 1200
\font\fivegoth=eufm5 scaled 1200
\newfam\gothfam
\textfont\gothfam=\tengoth \scriptfont\gothfam=\sevengoth
\scriptscriptfont\gothfam=\fivegoth

\def\ep{\varepsilon}

\begin{document}

\overfullrule=6pt
\title{
\bf
Asymptotics for multifactor Volterra type stochastic volatility models}
\author{
   { {\sc Giulia Catalini} } \and {\sc Barbara Pacchiarotti}\thanks{Dept. of Mathematics, University of Rome ``Tor Vergata''} }

\date{}
\maketitle

{\small
\bf Abstract. \rm We study multidimensional  stochastic volatility models in which the volatility process
is a positive continuous function  of a continuous multidimensional Volterra process that can be not self-similar. The
main results obtained in this paper are a generalization of the results due, in the one-dimensional case, to Cellupica and Pacchiarotti [M. Cellupica and B. Pacchiarotti (2021) Pathwise Asymptotics for Volterra Type Stochastic Volatility Models. {\em Journal of  Theoretical  Probability}, 34(2):682--727].
We state some (pathwise and finite-dimensional)
 large deviation  principles for the scaled log-price and as a consequence some
  (pathwise and finite-dimensional) short-time large deviation principles.}

\bigskip


\bigskip


\bigskip

{\small
\bf Corresponding Author\rm: Barbara Pacchiarotti,
Dipartimento di Matematica,
Universit\`a di Roma \sl Tor Vergata\rm, Via della Ricerca
Scientifica, I-00133 Roma, Italy.
}

\section{Introduction}\label{sect:intro}
The last few years have seen renewed interest in stochastic volatility
models  in which the volatility process
is a positive continuous function $\sigma$ of a continuous  stochastic process $\hat B$, that we assume to be  a Volterra type Gaussian process.
The goal of this paper is to extend to the multidimensional case and in a more general asset
a problem of large deviations for the log-price process of \textit{Volterra type stochastic volatility models},   studied in the one-dimensional case in \cite{Cel-Pac}, \cite{Gu2} and \cite{Gu3} (time-inhomogeneous case).
 Large deviations theory deals with the exponential decay of probabilities of \textquotedblleft rare events\textquotedblright, i.e. events whose probability is very small. These probabilities are important in many fields of study, including statistics, finance, engineering, statistical physics, and chemistry, since they often give informations about the large fluctuations of a random system around its most probable trajectory.
Jacquier, Pakkanen, and Stone \cite{Jac-Pak-Sto} prove a large deviations principle  (LDP)
for a scaled version of the log stock price process. In this same direction, Bayer, Friz, Gulisashvili,
Horvath and Stemper \cite{Bay-Fri-Gul-Hor-Ste}, Forde and Zhang \cite{For-Zha}, Horvath, Jacquier and Lacombe \cite{Hor-Jac-Lac}
and most recently Friz, Gassiat and Pigato \cite{Fri-Gas-Pig} (to name a few) prove large deviations principles
for a wide range of one-dimensional stochastic volatility models. Some results concerning asymptotic for the log-price processes,  in the multifactor setup can be found, for example, 
in \cite{Jac-Pan} where a sample path LDP for
log-processes associated with general Volterra systems is studied and in  \cite{Gu4} where a  comprehensive sample path LDP for
log-processes associated with multivariate time-inhomogeneous stochastic volatility models is stated. The main difference with the other models is that in this paper we can obtain short-time large deviations  even if  the Volterra process that appears in the variance is not self-similar (see Section \ref{sect:short}).
In the model considered in \cite{Cel-Pac}, the volatility is a positive continuous function of a continuous Volterra stochastic process. We generalize this model in a multi-factor setup, in which the volatility is a matrix of continuous function of a multidimensional continuous Volterra process, and we find a pathwise LDP for the family of the multivariate log-price process.
The generalization to the multidimensional case requires many technical details which as far as it seems to us cannot be cut out. So despite the many similarities with the previous paper, we are unable to write the proofs more concisely here.
Some recent literature uses stochastic volatility models with jumps (see, for example, \cite{Lin-Sen}, \cite{Sen-Wil-Nga},  \cite{Sal-Sen} and references therein), but as far as we know, there are no large  deviation results in these cases and this could be an interesting  topic for a future work.
In this paper, we consider a  quite general model in which the dynamic of the $\R^d$-valued 
process $(S(t))_{t\in[0,T]}$ $=(S_1(t),\ldots,S_d(t))_{t\in[0,T]}$ is modeled by the following equations:
\begin{equation}\label{eq:general-model}
	\begin{cases}
		\frac{dS_i(t)}{S_i(t)}=\mu_i(\hat{B}(t))\,dt+\sum_{\ell=1}^p\tilde{\sigma}_{i\ell}(\hat{B}(t))\,dB_\ell(t)+\sum_{j=1}^d\sigma_{ij}(\hat{B}(t))\,dW_j(t)\\
		S_i(0)=s_i^0
	\end{cases}
\end{equation}
for every $1\leq i	\leq d$, where $s^0\in\R^d$ is the initial value, $T>0$ is the time horizon and the process $\hat{B}$ is a non-degenerate continuous $\mathbb{R}^p$-valued multidimensional Volterra type process (see Definition \ref{def:Volterra-process}). 
The process $W$ is a $d$-dimensional standard Brownian motion independent from a $p$-dimensional standard Brownian motion $B$ appearing in definition of the process $\hat{B}$.
 This model has its own interest from a mathematical point of view due to its generality. Other large deviations results for multidimensional models can be found, for example,  in \cite{Chi-Fis} (for multidimensional diffusions) and \cite{Nua-Rov} (for multidimensional stochastic Volterra equations). Furthermore, the 
 processes $S_i$, $1\leq i	\leq d$, under suitable hypotheses on the coefficients can be  interpreted as price processes of correlated risky assets.
More precisely, the $d$ components of the process $S$ model the (dependent) prices of $d$ assets on the market. 
In this case the stochastic differential equation should be written as 
$$\frac{dS_i(t)}{S_i(t)}=\mu_i(\hat{B}(t))\,dt+\sum_{j=1}^d\Lambda_{ij}(\hat{B}(t))\,d\hat{W}_j(t)$$
with $\Lambda_{ij}$ some suitable functions, for every  $1\leq i,j\leq d$ and $\hat{W}$ a $d$-dimensional Brownian motion.
The matrix $\Lambda=(\Lambda_{ij})_{i,j=1,\ldots,d}$ is the volatility map and the matrix-valued process $(\Lambda(\hat B_t))_{t\in[0,T]}$ represents the joint volatility of the $d$ assets.
With a little abuse of language we call the process  $S$, defined in \eqref{eq:general-model}, the price process and $Z$ defined as
 $Z_i = \log S_i$  for $1 \leq i \leq d$ the log-price associated. 
 The initial condition for the log-process is denoted
by $x^0$. It is clear that $x^0_i=\log s^0_i$. 

It is assumed that $\mu_i:\mathbb{R}^p\to\mathbb{R}$,  $\sigma_{ij}:\mathbb{R}^p\to \mathbb{R}$ and $\tilde{\sigma}_{i\ell}:\mathbb{R}^p\to \mathbb{R}$ are continuous functions, for every $1\leq i,j\leq d$ and $1\leq \ell\leq p$.

The model is called uncorrelated when $\tilde{\sigma}_{i\ell}=0$, for every $1\leq i\leq d$ and $1\leq \ell\leq p$, i.e. when the dynamic of $S$ is not driven by $B$, otherwise it is called correlated. It is easy to understand what we mean with \textquotedblleft generalized model\textquotedblright; indeed, we find the one-dimensional model   taking $d=p=1$, $\tilde{\sigma}_{11}=\rho\sigma$ and $\sigma_{11}=\bar{\rho}\sigma$.
Notice that the one-dimensional case also generalize the model in \cite{Cel-Pac}.

We consider a suitable scaled version $(Z^n)_n$ (see equation \eqref{eq:scaled-log-price-corr}) of the log-price process $Z$ and
we  obtain a  sample path  large deviation  principle  for the family  of  processes $ ((Z_t^{n}-x_0)_{t \in [0,T]})_{n \in \mathbb{N}}.$
A large deviation  principle  for $ (Z_T^{n}-x_0)_{n \in \mathbb{N}} $ can be  obtained with the same techniques, but  it can be also obtained
by contraction and this is the approach we follow here. We always suppose without loss of generality that $s_i^0=1$ (and then $x_i^0=0$), for every $1\leq i\leq d$.

The paper is organized as follows.
In Section \ref{sect:ldp} we recall some basic facts about large deviations for continuous Gaussian processes  and we give the definition of multidimensional Volterra process.
For  some facts about large deviations   for joint and marginal distributions we refer for simplicity to Section 3 in \cite{Cel-Pac} that is   a summary of the results we use and that are proved in \cite{Cha}.
In Sections \ref{sect:uncorrelated} and \ref{sect:correlated} are contained the main results.
More precisely in Section \ref{sect:uncorrelated} we prove a large deviation  principle for the log-price process in  the uncorrelated model.
In Section \ref{sect:correlated}, we prove a large deviation  principle for the log-price process in  the correlated model 
(we will follow the same pattern as in  \cite{Cel-Pac}). In Section  \ref{sect:short}, following \cite{GiPaPi} we obtain a multidimensional short-time LDP.

\section{LDP for multidimensional Volterra processes} \label{sect:ldp}

We briefly recall some main facts on large deviation  principles and Volterra processes   we are going to use.
For a detailed development of this very wide theory we can refer, for example, to the following classical references:
the book of Varadhan \cite{Var},
Chapitre II in Azencott \cite{Aze},
Section 3.4 in Deuschel and Strook \cite{Deu-Str}, Chapter 4 (in particular Sections 4.1, 4.2  and 4.5) in  Dembo and Zeitouni
\cite{Dem-Zei}, for large deviation  principles;   \cite{De-Ust} and  \cite{Hult} for
Volterra processes.

\subsection{Large deviations}

\begin{definition}
Let $E$ be a topological space,   $\cl{B}(E)$ the Borel $\sigma$-algebra and  $(\eta_n)_{n\in \N}$ a family of probability measures on $\cl{B}(E)$; let $\gamma \, : \N \rightarrow \R^+ $ be a speed function, i.e. $\gamma_n \rightarrow +\infty$ as $n\to +\infty$. We say that the family of probability measures $(\eta_n)_{n\in \N}$ satisfies a LDP  on $E$ with the rate function $I$ and the speed $\gamma_n$ if, for any open set
$\Theta$,
$$
-\inf_{x \in {\Theta} } I(x) \le \liminf_{n\to+\infty} \frac 1{\gamma_n} \log \eta_n (\Theta)
$$ and for any closed set $\Gamma$
\begin{equation} \label{eq:upperbound} \limsup_{n\to +\infty}\frac 1 {\gamma_n} \log \eta_n (\Gamma) \le -\inf_{x \in {\Gamma}} I(x).
\end{equation}
\end{definition}
A rate function is a lower semicontinuous mapping $I:E\rightarrow [0,+\infty]$. A rate function $I$ is said \textit{good} if  $\{I\le a\}$ is a compact set for every $a \ge 0$.
\begin{definition}
Let $E$ be a topological space,   $\cl{B}(E)$ the Borel $\sigma$-algebra and  $(\eta_n)_{n\in \N}$ a family of probability measures on $\cl{B}(E)$; let $\gamma \, : \N \rightarrow \R^+ $ be a speed function. We say that the family of probability measures $(\eta_n)_{n\in \N}$ satisfies a weak LDP (WLDP) on $E$ with the rate function $I$ and the speed $\gamma_n$ if the upper bound (\ref{eq:upperbound}) holds for compact sets.
\end{definition}
From now on, given $ T>0 $, we will denote with $ \cl C ^p $ the set of $\R^p$-valued continuous functions on $ [0,T] $,  $ \cl C  $ if $p=1$, endowed with the topology induced by the sup-norm $\|\cdot\|_\infty$, i.e if $f=(f_1,\ldots,f_p)$ then
$$\|f\|_\infty=\sup_{0\leq t \leq T}||f(t)||,$$
where $||\cdot||$ is the euclidean norm in $\R^p$.
We will denote with $  \cl C_0 ^p $, $ \cl C_0  $ if $p=1$, the subspace of $\R^p$-valued continuous functions on $ [0,T] $ starting from zero.
In what follows, we will always suppose our processes to be continuous.
\begin{remark}\rm
We say that a family of continuous processes $((U^n(t))_{t \in [0,T]})_{n\in \N}$, $U^n(0)=0$ satisfies a LDP if  the family of their laws satisfies a LDP on $  \cl C ^p$.
\end{remark}

Let us conclude this section with some important definitions.
\begin{definition}
\label{def:exp-equiv}
Let $ (E, d_{E}) $ be a metric space (we consider on $E$ the Borel $\sigma$-algebra) and let $ (Z^n)_{n\in\N} $ and $ (\tilde{Z}^n)_{n\in\N} $ be two families of $ E $-valued random variables.  Then $ (Z^n)_{n\in\N} $ and $ (\tilde{Z}^n)_{n\in\N} $ are exponentially equivalent (at the speed $\gamma_n$) if  for any $  \delta > 0 $,
	$$ \limsup_{n\to +\infty }\frac{1}{\gamma_n} \log P(d_E(\tilde{Z}^n,Z^n)>\delta) = -\infty. $$
	
\end{definition}
As far as the LDP is concerned exponentially equivalent measures are indistinguishable. See  Theorem 4.2.13 in \cite{Dem-Zei}.

 Let us give the definition of exponentially good approximation and some results about this topic. Here the main reference  is \cite{Bax-Jai}.

\begin{definition}\label{def:exp-approx}
	Let $ (E,d_E) $ be a metric space; consider  the $ E $-valued random variables $ Z^n $ and $ Z^{n,m}$. The families $ (Z^{n,m} )_{n\in \N}$ for $m\geq 1$ are called   \textbf{exponentially good approximations} of $( {Z}^{n})_{n\in \N} $ at the speed $\gamma_n$ if, for every $ \delta>0 $,
	$$  \lim_{m \to +\infty} \limsup_{n\to +\infty}\frac{1}{\gamma_n} \log \mathbb{P}(d_E(Z^{n,m}, Z^n) > \delta)= -\infty. $$
	
\end{definition}

Next theorem, Theorem 3.11 in \cite{Bax-Jai}, states that under a suitable condition if for each $ m\geq 1 $ the sequence $ (Z^{n,m})_{n\in \N} $ satisfies a LDP with the rate function $ I^m ,$ then also $ (Z^{n})_{n\in \N} $  satisfies a LDP with the rate function $ I$,  obtained in terms of the $ I^m(\cdot) $'s.

\begin{theorem}\label{th:exp-good-approx}{\rm[Theorem 3.11 in \cite{Bax-Jai}]}
Let  $ (Z^{n,m} )_{n\in \N}$ for $m\geq 1$ be  exponentially good approximations of $( {Z}^{n})_{n\in \N} $ at the speed $\gamma_n$. Suppose that $ (Z^{n,m})_{n\in \N} $ satisfies a LDP with the speed $ \gamma_n $ and the good rate function $ I^m$.
	Then $(Z^n)_{n\in \N} $ satisfies a LDP with the speed $ \gamma_n $ and the good rate function $ I $ given by
	$$ {I}(x)=\lim_{\delta \to 0}\liminf_{m \to +\infty}\inf_{y \in B_\delta(x)}I^m(y)=\lim_{\delta \to 0}\limsup_{m \to +\infty}\inf_{y \in B_\delta(x)}I^m(y),$$
with $  B_\delta(x)=\{y \in E: d_E(x,y)<\delta\} $
	\end{theorem}
We also need the following proposition.
\begin{proposition}\label{prop:identification}{\rm[Proposition 3.16 in \cite{Bax-Jai}]}
	In the same hypotheses of   Theorem \ref{th:exp-good-approx}, if
\begin{itemize}
\item $ I^m(x)\overset {m \rightarrow +\infty}{{\longrightarrow}} J(x)$, for   $x\in E$;
\item $ x_m \overset{m \rightarrow +\infty}{{\longrightarrow}}x$ implies $  \liminf_{m \to+\infty}I^m(x_m)\geq J(x),$
\end{itemize}
for some functional $J(\cdot)$, then $ I(\cdot)=J(\cdot)$.
\end{proposition}

\subsection{Multidimensional Volterra  processes}

Let $ (\Omega,\cl{F},(\cl{F}_t)_{t\in[0,T]},\mathbb{P}) $ be a probability space and $ B=(B(t))_{t \in [0,T]} $ a standard Brownian motion.
Suppose $ \hat{B}=(\hat{B}(t))_{t \in [0,T]} $ is a centered (real) Gaussian process having the following Fredholm representation,
\begin{eqnarray}\label{eq:integral-representation} \hat{B}(t)= \int_0^T K(t,s)\, dB(s), \quad 0 \leq t \leq T,
\end{eqnarray}
where  $ T > 0 $ and  $ K $ is a  measurable square integrable kernel on $ [0,T]^2 $ such that
$$\sup_{t\in[0,T]} \int_0^T K(t,s)^2 \,ds < \infty. $$

The modulus of continuity of the kernel $ K $  is defined as follows
$$
M(\delta)=  \sup_{\{t_1,t_2 \in [0,T]: |t_1-t_2|\leq \delta\}}  \int_0^T |K(t_1,s)-K(t_2,s)|^2\,ds, \quad 0 \leq \delta \leq T.
$$

  The covariance function of the process $ \hat{B} $ is given by
$$k(t,s)= \int_0^T K(t,u)K(s,u)\,du, \quad t,s \in [0,T].$$
Let us define a Volterra process.	
\begin{definition}\label{def:Volterra-process}
	The process in (\ref{eq:integral-representation}) is called a Volterra type Gaussian process if the following conditions hold for the kernel $ K $:
	\begin{enumerate}
		\item [\textit{(a)}] $K(0,0)=0$  and $ K(t,s)=0 $ for all $ 0\leq t < s \leq T $;
		\item [\textit{(b)}] There exist constants $ c > 0 $ and $ \alpha > 0 $ such that $ M(\delta)\leq c\,\delta^\alpha $ for all $ \delta \in [0,T] $.
		\end{enumerate}
\end{definition}

\begin{remark}\rm
Condition \textit{(a)} is a typical Volterra type condition for the kernel $ K $ and the integral representation in (\ref{eq:integral-representation}) becomes $\hat{B}(t)= \int_0^t K(t,s)\, dB(s),$ for $ 0 \leq t \leq T. $ So $ \hat{B} $ is adapted to the natural filtration generated by $B$.
		Condition \textit{(b)} guarantees the existence of a H\"{o}lder continuous modification of the process $ \hat{B} $.
\end{remark}

In this section we introduce $\R^p$-valued Gaussian processes with independent components in which each component is a Volterra type Gaussian process.
\begin{definition}
	The process $\hat{B}=(\hat{B}(t))_{t\in [0,T]} $ defined by:
	\begin{eqnarray}\label{eq:volterra-multidimensional}
		\hat{B}(t)=(\hat{B}_1(t),\ldots,\hat{B}_p(t))=\Big(\int_{0}^{T}K_1(t,s)\, dB_1(s),\ldots,\int_{0}^{T}K_p(t,s)\, dB_p(s)\Big)
	\end{eqnarray}
	for every  $0 \leq t \leq T$, is a \textbf{multidimensional Volterra type process}
	if for every $1\leq i\leq p$,  $K_i$ satisfies conditions (a), (b) in Definition \ref{def:Volterra-process}.
\end{definition}

\begin{remark}\rm
	The process $\hat{B}$ defined in (\ref{eq:volterra-multidimensional}) is a centered Gaussian process with covariance function $k:[0,T]\times [0,T]\to\R^{p\times p}$ given by
	\begin{eqnarray*}
		k(t,s)_{ij}=\delta_{ij} \int_{0}^{t \wedge s}K_i(t,u)K_j(s,u)\, du
	\end{eqnarray*}
	for every $1\leq i,j\leq p$ and $t,s \in [0,T]$.
Moreover, $\hat{B}$ admits a continuous version with H\"{o}lder continuous sample paths of index $\gamma$ for every $\gamma < \alpha/2$.
\end{remark}

Let $(\hat{B}^n)_{n\in\N}$ a family of processes such that, for every $n\in\N$, $(\hat{B}^n(t))_{t\in [0,T]}$ is a continuous Volterra type  process of the form
$$
		\hat{B}^n(t)=(\hat{B}^n_1(t),\ldots,\hat{B}^n_p(t))=\Big(\int_{0}^{T}K_1^n(t,s)\, dB_1(s),\ldots,\int_{0}^{T}K_p^n(t,s)\, dB_p(s)\Big)
	$$
	with $K_i^n$  suitable kernels and covariance function
	$$k^n_{ij}(t,s)=\delta_{ij}\int_0^{t\wedge s}K^n_i(t,u)K^n_j(s,u)\, du\quad s,t\in[0,T], \quad i,j=1,\ldots,p.$$

Suppose that there exist   kernels $ K_\ell $ such that, for every $\ell=1,\ldots,p$,
\begin{equation} \label{eq:ker-limit} \lim_{n \to +\infty} \frac{K^n_{\ell}(t,s)}{{\ep_n}}=K_{\ell}(t,s). \end{equation} Under
suitable   conditions on the covariance functions, a  LDP for $((\ep_n B, \hat B^n))_{n\in \N}$ holds.  More precisely,
 $ ((\ep_nB, \hat{B}^n))_{n \in \mathbb{N}} $ satisfies a LDP on $( \cl C_0 ^p)^2$ with the speed $ \ep_n^{-2} $ and the good rate function
	\begin{eqnarray}\label{eq:rate-function-B-hatB}
		I_{(B,\hat{B})} (f,g)=
		 \begin{cases}
		  \frac12  \int_0^T ||\dot{f}(s)||^2 \, ds & (f,g) \in \cl H^p_{(B,\hat B)}\\
		  +\infty & otherwise
		\end{cases}
		\end{eqnarray}
where
\begin{equation}\label{eq:RKHS-B-hatB}
	\mathscr{H}_{(B,\hat{B})}^p=\Big\{ (f,g)\in ( \cl C_0 ^p)^2 : f\in H_0^{1,p}[0,T], g=\hat{f}\Big\},
\end{equation}
\begin{equation}\label{eq:hatf}
\hat f=(\hat f_1,\ldots, \hat f_p), \quad  \hat f_\ell(t)=\int_0^t K_\ell(t,u)\dot f_\ell (u) \, du \quad t\in[0,T],\,  1\leq \ell \leq p,\end{equation}
$K_\ell$ is defined from equation (\ref{eq:ker-limit}) and $H_0^{1,p}[0,T]$ is the Cameron-Martin space of the $p$-dimensional Brownian motion.

In this paper we are not interested to this problem and we assume that such LDP holds. For details about this LDP for Gaussian processes see \cite{Mac-Pac}, \cite{Gio-Pac} and \cite{Cel-Pac}.
\begin{assumption}\label{ass:LDP-B-hatB}
$ ((\ep_nB, \hat{B}^n))_{n \in \mathbb{N}} $ satisfies a LDP on $( \cl C_0 ^p)^2$ with the speed $ \ep_n^{-2} $ and the good rate function
given by equation (\ref{eq:rate-function-B-hatB}).
	\end{assumption}
	
\begin{example}\rm
Let us consider two examples in which Assumption \ref{ass:LDP-B-hatB} is satisfied.
 We consider a multidimensional version of the models studied in  \cite{GiPaPi}.
 
 Let us define the multidimensional log fractional Brownian motion as 
 the process defined in \eqref{eq:volterra-multidimensional} 
	with $$K_\ell(t,s)=C (t-s)^{H-1/2} (-\log(t-s))^{-a},\quad \ell=1,\ldots,p$$
for $0\leq H \leq 1/2$,  $a>1$ and $C>0$ a positive constant.
	
 For $\eta_n\to 0$ consider $\hat{B}^n(t)=\hat{B}(\eta_n t)$. Then
 $$
\hat{B}^n(t)=\hat{B}(\eta_n t)=\Big(\int_{0}^{
 t}K^n( t,s)\, dB_1(s),\ldots,\int_{0}^{ t}K^n
( t,s)\, dB_p(s)\Big),
$$
where $$K^n(t,s)=\sqrt{\eta_n} K(\eta_n t,\eta_n s).$$
Then	$ ((\ep_nB, \hat{B}^n))_{n \in \mathbb{N}} $ satisfies a LDP on $( \cl C_0 ^p)^2$ with the speed $ \ep_n^{-2}= \eta_n^{-2H}(-\log \eta_n)^{2p}$ and the good rate function defined in \eqref{eq:ker-limit} with \it limit kernel \rm that is the one of the Riemann-Liouville fractional Brownian motion.

 In the same way let us define the multidimensional fractional  Ornstein-Uhlenbeck process.
	Here $$K_\ell(t,s)=K_H(t,s)- a \int_s^t e^{-a(t-u)} K_H(u,s) du, \quad \ell=1,\ldots,p$$
for $a>0$,  $H\in(0,1)$ and $K_H$ the kernel of the fractional Brownian motion.
	
 For $\eta_n\to 0$ consider
 $
\hat{B}^n(t)=\hat{B}(\eta_n t)$.
Then	$ ((\ep_nB, \hat{B}^n))_{n \in \mathbb{N}} $ satisfies a LDP on $( \cl C_0 ^p)^2$ with the speed $ \ep_n^{-2}= \eta_n^{-2H}$ and the good rate function defined in \eqref{eq:ker-limit} with \it limit kernel \rm that is the one of the  fractional Brownian motion.
\end{example}

\begin{remark}\label{rem:LDP-hatB}
From Assumption \ref{ass:LDP-B-hatB} and the contraction principle, the family $ (\hat{B}^n)_{n \in \mathbb{N}} $ satisfies a LDP on $  \cl C_0 ^p $
with the speed $ \ep_n^{-2} $ and the good rate function
\begin{eqnarray}\label{eq:rate-function-hatB}
I_{\hat{B}}(g)=\inf\Big\{\frac12 \int_0^T ||\dot{f}(s)||^2 \, ds  : \hat f=g, \,\, f\in H_0^{1,p}[0,T]\Big\},
\end{eqnarray}
with the understanding $ I_{\hat{B}}(g)= +\infty $ if the set
is empty.
\end{remark}

	\section{The uncorrelated model} \label{sect:uncorrelated}
	The results of this section are an intermediate step in order to  obtain the general case. Here it is assumed that the prices of the $d$ assets are independent.
The dynamic of in the uncorrelated model of  $(S(t))_{t\in[0,T]}=(S_1(t),\ldots,S_d(t))_{t\in[0,T]}$ is the one  in \eqref{eq:general-model} with $\tilde\sigma_{i\ell}=0$ for every $1\leq i\leq d$, $1\leq \ell \leq p$, that is
\begin{equation}\label{eq:asset-prices-uncorr}
		\begin{cases}
			dS_i(t)=S_i(t)\mu_i(\hat{B}(t))\,dt+S_i(t)\sum_{j=1}^{d}\sigma_{ij}(\hat{B}(t))\,dW_j(t), \quad 0\leq t\leq T \\ S_i(0)=s_i^0
		\end{cases}
	\end{equation}
for every $1\leq i	\leq d$, where, we recall, $s^0=(s^0_1,\ldots,s^0_d)\in\R^d$ is the initial price, $T>0$ is the time horizon and the process $\hat{B}$ is a non-degenerate continuous multidimensional Volterra process.
The unique solution to the equation   is the  exponential
$$S_i(t)=s_i^0\exp\Big\{\int_{0}^t\Big(\mu_i(\hat{B}(s))-\frac12\sum_{j=1}^d\sigma_{ij}(\hat{B}(s))^2\Bigr)\, ds + \int_0^t\sum_{j=1}^d\sigma_{ij}(\hat{B}(s))\, dW_j(s)\Bigr\}$$
for every $1\leq i\leq d$ and $0\leq t \leq T$ (for further details see Section IX-2 in \cite{Re-Yor}).
Therefore, the log-price processes $X_i(t)=\log S_i(t)$, with $X_i(0)=x_i^0=\log s_i^0$ is
\begin{equation*}
	X_i(t)=x_i^0+ \int_0^t \Big(\mu_i(\hat{B}(s))-\frac12\sum_{j=1}^d\sigma_{ij}(\hat{B}(s))^2\Bigr)\, ds + \int_0^t \sum_{j=1}^d\sigma_{ij}(\hat{B}(s))\,dW_j(s),
\end{equation*}
for every $1\leq i\leq d$ and $0\leq t\leq T$.
For the sake of simplicity (it is not restrictive), from now on we assume that the initial conditions $s_i^0$ for the asset prices satisfy $s_i^0=1$ and so $x_i^0=\log s_i^0=0$ for every $1\leq i\leq d$.
 Now, let $\ep: \N\to\R_+$ be an infinitesimal  function i.e $\ep_n\to 0$, as $n\to+\infty$. For every $n\in\N$, we will consider the following scaled version of the stochastic differential equations in (\ref{eq:asset-prices-uncorr})
$$\begin{cases}
	dS_i^n(t)=S_i^n(t)\mu_i(\hat{B}^n(t))\,dt+\ep_nS_i^n(t)\sum_{j=1}^d\sigma_{ij}(\hat{B}^n(t))\,dW_j(t)\quad 0\leq t\leq T\\ S_i^n(0)=1
\end{cases}$$
for every $1\leq i\leq d$, where $(\hat{B}^n)_{n\in\N}$
is a family of multidimensional Volterra processes such that   Assumption \ref{ass:LDP-B-hatB} is fulfilled.
The log-price processes $X_i^n(t)=\log S_i^n(t)$ in the scaled model are given by
\begin{equation}\label{eq:scaled-log-price-uncorr}
	X_i^n(t)=\int_0^t\Big(\mu_i(\hat{B}^n(s))-\frac12\ep_n^2\sum_{j=1}^d\sigma_{ij}(\hat{B}^n(s))^2\Big)\, ds+\ep_n\int_0^t\sum_{j=1}^d\sigma_{ij}(\hat{B}^n(s))\,dW_j(s)
\end{equation}
for every $1\leq i\leq d$ and $0\leq t\leq T$.

We made the following assumptions on the coefficients.
\begin{assumption}\label{ass:hp-sigma-mu-I} 
 $\mu_i:\R^p\to\R$ and $\sigma_{ij}:\R^p\to \R$  are continuous functions, for every $1\leq i,j\leq d$ and
  $det(a(y))\neq0$ for every $y\in\R^p$ where 
 $a=\sigma\sigma^T$.
\end{assumption}

Note that the hypotheses of continuity on the coefficients $\sigma_{ij}$ are quite mild. Similar hypotheses can be found in \cite{Gu4} and,  for example, in  the classical
Bergomi model  (see \cite{Bay-Fri-Gat}), where the volatility is the exponential (i.e. a continuous function) of a Volterra process. The hypothesis on the determinant are the generalization to the multidimesional case of the request on the volatility map to be positive.
\begin{remark}\label{rem:unif-def-pos}\rm
	Under Assumption \ref{ass:hp-sigma-mu-I},  there exists the inverse matrix $a^{-1}(y)$ for every $y\in \R^p$ and $a^{-1}_{ij}$ are continuous functions for every $i,j=1,\ldots,d$.
Furthermore the matrix $a$ is uniformly strictly positive definite on compact sets, i.e.
for every unit vector  $x\in\R^d$ and $K\subset\R^p$ compact there exists $\alpha_K>0$ such that
$\inf_{y\in K}x^Ta(y)x\geq\alpha_K>0$, where $\alpha_K=\inf_{y\in K} \lambda_{min}(y)$, being $\lambda_{min}(y)$ the minimum eigenvalue of $a(y)$. Therefore the matrix $a(y)-\alpha_K I$ is
uniformly  positive definite on $K$.	And in a similar way we have that the matrix $\beta_K I-a(y)$, where
$\beta_K=\sup_{y\in K} \lambda_{max}(y)$, being $\lambda_{max}(y)$ the maximum eigenvalue of $a(y)$, is
uniformly  positive definite on $K$.  Note that under Assumption \ref{ass:hp-sigma-mu-I}   we have the same properties for the matrix $a^{-1}$.
\end{remark}
The main result of this section is the following 
 sample path LDP for the process $ (X^{n})_{n \in \N} $.  
\begin{theorem}\label{th:main-uncorr}
Under  Assumptions \ref{ass:hp-sigma-mu-I} and \ref{ass:LDP-B-hatB}
	a  LDP with the speed $ \ep^{-2}_n $ and the good rate function
$ I_X(\cdot)$ defined in (\ref{eq:log-price-rate-function-uncorr})
	holds for the family $ (X^{n})_{n \in \N} $ of processes defined in (\ref{eq:scaled-log-price-uncorr}).
\end{theorem}

In order to prove this theorem we procede as in \cite{Cel-Pac}. Checking that hypotheses of  Theorem 3.3.(Chaganty's Theorem) in   \cite{Cel-Pac}) are fullfilled we obtain a WLDP for the family $(\hat B^n, X^n)_n\in \N$ with a rate function $I(\,\cdot\,|\,\cdot\,)$; then we prove that $I(\,\cdot\,|\,\cdot\,)$ is a good rate function and therefore (by 
Chaganty main result) the family $(X^n)_n\in \N$ satisfies a (full) LDP. 
Assumption \ref{ass:LDP-B-hatB} ensures that condition (i) in Theorem 3.3 is fullfilled more precisely, 
 the family $(\hat B^n)_n\in \N$ satisfies a LDP on $  \cl C_0 ^p $ with  the  speed $ \ep_n^{-2}$ and the good rate function $ I_{\hat{B}}(\cdot) $ given by (\ref{eq:rate-function-hatB}). So we have to prove only condition (ii) in Theorem 3.3, that  is the LDP continuity condition (see Definition 3.2 in \cite{Cel-Pac}).
  For the multidimensional case we   need some technical results  on positive definite matrices which we postpone in the appendix.

 We want to investigate the behavior of $X^n$, when conditioned to the process $\hat B^n$. More precisely, we want to establish a LDP for the family of the conditional processes
$$X^{n,\varphi}=X^n|\big(\hat{B}^n(t)=\varphi(t)\quad 0\leq t\leq T\big)$$
as $n\to+\infty$.
For (almost) every $ \varphi \in  \cl C_0 ^p $,  we have (in law)
\begin{equation}\label{eq:logprice-uncorr-scaled-cond}
X_i^{n,\varphi}(t)=\int_0^t\Big(\mu_i(\varphi(s))-\frac12\ep_n^2\sum_{j=1}^d\sigma_{ij}(\varphi(s))^2\Big)\, ds + \ep_n\int_0^t\sum_{j=1}^d\sigma_{ij}(\varphi(s))\, dW_j(s)\end{equation}
for every $1\leq i\leq d$ and $0\leq t\leq T$.
It is enough to show that the conditions $(a)$, $(b)$ and $(c)$ of LDP continuity condition are satisfied (condition (ii) of Chaganty's Theorem).
The proofs of the next results are quite similar to the proofs of Proposition  5.4, Proposition 5.6 and Lemma 5.7 in \cite{Cel-Pac}. Anyway  there are many technical issues related to the multidimensional case therefore we give here all the details.
Let us define the following functional $\Gamma$  that will be useful in the sequel.

For     $A\in {\cl C}^{d\times d}$ symmetric and positive definite (i.e. $A(t)$ is a symmetric and positive definite matrix for every $t\in[0,T]$) and
 $x\in  \cl C^d$   define the functional
$$\Gamma(x|A)=\begin {cases}\frac12\int_0^T\dot{x}(t)^TA(t)\dot{x}(t)\,dt&x\in  H_0^{1,d}[0,T]\\
+\infty & otherwise\end{cases}$$

\begin{remark}\label{rem:Gamma}
For  $x,y,z\in  H_0^{1,d}[0,T]$ and $A, B\in {\cl C}^{d\times d}$ we have,
	\begin{enumerate}
		\item[\textit{(i)}] $\Gamma(x|B)\geq\Gamma(x|A)$ if $B-A$ is  positive definite;
		\item[\textit{(ii)}]
$\Gamma(x+y|A)\leq 2\Gamma(x|A)+2\Gamma(y|A)$;
\item[\textit{(iii)}]
		$\Gamma(x+y+z|A)\leq 3\Gamma(x|A)+3\Gamma(y|A)+3\Gamma(z|A). $
\end{enumerate}
\end{remark}

\begin{proposition} \label{prop:LDPcc-uncorr}
Under  Assumptions \ref{ass:hp-sigma-mu-I}
the sequence of the conditional processes $ (X^{n,\varphi})_{n \in \N} $  satisfies the LDP continuity condition  with the  speed $ \ep_n^{-2}$ and the good rate function
\begin{equation}\label{eq:cond-rate-function-uncorr}\begin{split}
		  J(x|\varphi)&=\Gamma\Big(x-\int_0^\cdot\mu(\varphi(t))\,dt\Big| a^{-1}(\varphi)\Big)\\&=\begin{cases} \frac12\int_0^T(\dot{x}(t)-\mu(\varphi(t)))^Ta^{-1}(\varphi(t))(\dot{x}(t)-\mu(\varphi(t)))\,dt&x\in H_0^{1,d}[0,T]\\
			+\infty&otherwise.
		\end{cases}\end{split}\end{equation}
 \end{proposition}

\proof
$(a)$ 	
For every $\varphi\in  \cl C_0 ^p$  we  prove that $ (X^{n,\varphi})_{n \in \N} $ obeys a LDP  by using the multidimensional version
of Theorem 2.14 in \cite{Cel-Pac} that still holds as a simple consequence of Theorem 3.1 in \cite{Chi-Fis}. For every $n\in\N$ and $t\in [0,T]$, the coefficients are:
	\begin{itemize}
		\item[-]$b_n(t)=\mu(\varphi(t))-\displaystyle\frac12\ep_n^2 Diag(a(\varphi(t)))$ (with $Diag(C)$ we mean the diagonal of the matrix $C$) and so $b_n\to \mu(\varphi)$ in $ \cl C^d$, as $n\to+\infty$.
		\item[-]$\sigma_n(t)=\sigma(\varphi(t))$, for every $t\in [0,T]$, not depending on $n$.
		\end{itemize}
	
Then the family $ (X^{n,\varphi})_{n \in \N} $ satisfies a  LDP  with the  speed $ \ep^{-2}_n $ and the good rate function
$$J(x|\varphi) = \inf \Big\{\frac 12  \int_0^T ||\dot{y}(t)||^2 \, dt:\, \int_0^t\mu(\varphi(s))\,ds+\int_0^t\sigma(\varphi(s))\dot{y}(s)\, ds=x(t), \,  y \in H_0^{1,d}[0,T]\Big\}$$
with the usual  understanding $J(x|\varphi)= +\infty  $ if the set is empty.
If $ y \in H_0^{1,d}[0,T]$ then
$$ \dot{x}(t)=\mu(\varphi(t))+\sigma(\varphi(t))\dot{y}(t) \,\,\mbox{ a.e., with } x(0)=0.$$
Thanks to Remark \ref{rem:continuous} $(ii)$ and Assumption \ref{ass:hp-sigma-mu-I} we have that  $det(\sigma\circ \varphi)\neq 0$ and the rate function above simplifies to $ J( \cdot|\varphi) $ defined in (\ref{eq:cond-rate-function-uncorr}).
\medskip	

$(b)$
	Let $(\varphi_n)_{n\in\N}\subset  \cl C_0 ^p$ such that $\varphi_n\to \varphi$ in $ \cl C_0 ^p$, as $n\to+\infty$. For every $n\in\N$ we consider the law of the $\R^d$-valued Gaussian diffusion process $(X^{n,\varphi_n}(t))_{t\in[0,T]}$ defined by (\ref{eq:logprice-uncorr-scaled-cond}), where the function $\varphi$ is replaced  by  $\varphi_n$.
	 The coefficients now are:
		\begin{itemize}
		\item[-]$b_n(t)=\mu(\varphi_n(t))-\displaystyle\frac12\ep_n^2 Diag(a(\varphi_n(t)))$ , for every $t\in [0,T]$;
		\item[-]$\sigma_n(t)=\sigma(\varphi_n(t))$, for every $t\in [0,T]$.
	\end{itemize}
From Remark \ref{rem:conv-Mphi-Aphi-AphiInv}, $a_{ii}(\varphi_n)\to a_{ii}(\varphi)$ in $\cl C$, as $n\to+\infty$, for every $1\leq i\leq d$. Therefore $b_n(t)\to\mu(\varphi(t))$ and $\sigma_n(t)\to\sigma(\varphi(t))$ so $b_n\to \mu(\varphi)$ in $ \cl C^d$ and  $\sigma_n\to \sigma(\varphi)$ in ${\cl C}^{d\times d}$,
as $n\to+\infty$. Thus
 $(X^{n,\varphi_n})_{n\in\N}$ obeys a LDP with the  speed $ \ep^{-2}_n $ and the good rate function $ J(\cdot|\varphi)$ defined in
(\ref{eq:cond-rate-function-uncorr}).

\medskip

$(c)$ Now, we will check that $J(\,\cdot\,|\,\cdot\,)$ is lower semi-continuous as a function of $(\varphi,x)\in  \cl C_0 ^p\times { \cl C_0 }^d$.
Let $(\varphi_n,x_n)_{n\in\N}\subset  \cl C_0 ^p\times { \cl C_0 }^d$ and $(\varphi,x)\in  \cl C_0 ^p\times { \cl C_0 }^d$ be functions such that
	$(\varphi_n,x_n)\to(\varphi,x)$
	in $ \cl C_0 ^p\times { \cl C_0 }^d$, as $n\to+\infty$.
	If $$\displaystyle\liminf_{n\to +\infty}J(x_n|\varphi_n)=\lim_{n\to +\infty}J(x_n|\varphi_n)=+\infty,$$ there is nothing to prove. Therefore, up to a subsequence, we can suppose that $(x_n)_{n\in\N}\subset H_0^{1,d}[0,T]$.
	 Now, from Lemma \ref{lemma:tech2} and the definition of $\Gamma$,  we have that for every $\varepsilon>0$ there exists $n_\varepsilon \in\N$ such that for every $n\geq n_\varepsilon$
	\begin{eqnarray*}
		J(x_n|\varphi_n)&=&	\Gamma\Big(x_n-\int_0^\cdot\mu(\varphi_n(t))\,dt\Big| a^{-1}(\varphi_n)\Big)
		\geq (1-\varepsilon)	\Gamma\Big(x_n-\int_0^\cdot\mu(\varphi_n(t))\,dt\Big| a^{-1}(\varphi)\Big)\\
		&=&(1-\varepsilon)J\Big(x_n-\int_0^\cdot\mu(\varphi_n(t))\,dt+\int_0^\cdot\mu(\varphi(t))\,dt\Big| \varphi \Big).
	\end{eqnarray*}
	Moreover thanks to Remark \ref{rem:conv-Mphi-Aphi-AphiInv},  we have
	$$x_n-\int_0^\cdot\mu(\varphi_n(t))\,dt+\int_0^\cdot\mu(\varphi(t))\,dt\to x$$
	in ${ \cl C_0 }^d$, as $n\to+\infty$. The functional $J(\,\cdot\,|\varphi)$, being a rate function, is lower semi-continuous on ${ \cl C_0 }^d$, therefore
	\begin{eqnarray*}
		\liminf_{n\to +\infty}J(x_n|\varphi_n)
		&\geq&(1-\varepsilon)\liminf_{n\to +\infty} J\Big(x_n-\int_0^\cdot\mu(\varphi_n(t))\,dt+\int_0^\cdot\mu(\varphi(t))\,dt\big|\varphi\Big)\geq (1-\varepsilon) J(x|\varphi)
	\end{eqnarray*}
for every $\varepsilon>0$ and thus we can conclude.
\cvd

\begin{proposition} Suppose that Assumption \ref{ass:LDP-B-hatB} holds and
 that $\sigma$ and $\mu$ satisfy Assumption \ref{ass:hp-sigma-mu-I}.
	 Then $(\hat{B}^n,X^n)_{n\in\N}$ satisfies the WLDP with the rate function
\begin{equation}\label{eq:rate-function-couple}
		I(\varphi,x)=I_{\hat{B}}(\varphi)+J(x|\varphi)
\end{equation}
for $(\varphi,x)\in  \cl C_0 ^p\times { \cl C_0 }^d$, where $I_{\hat{B}}(\,\cdot\,)$ and $J(\,\cdot\,|\,\cdot\,)$ are defined in  (\ref{eq:rate-function-hatB}) and (\ref{eq:cond-rate-function-uncorr}), respectively.
	 Moreover, $(X^n)_{n\in\N}$ satisfies the LDP  with the  speed $\ep_n^{-2}$ and the rate function
	\begin{equation}\label{eq:log-price-rate-function-uncorr}
		I_X(x)=\begin{cases}\displaystyle\inf_{f\in H_0^{1,p}[0,T]}\Big\{ \frac12 \|f\|_{H_0^{1,p}[0,T]}^2+J(x|\hat{f}) \Big\}\quad &x\in H_0^{1,d}[0,T]\\
			\phantom{ff}+\infty&otherwise
		\end{cases}
	\end{equation}
where $\hat{f}$ is defined in (\ref{eq:hatf}).
\end{proposition}

\proof
Thanks to Remark \ref{rem:LDP-hatB}, the family  $(\hat B^n)_n$  satisfies a LDP on $  \cl C_0 ^p $
with the speed $ \ep_n^{-2} $ and the good rate function
\begin{eqnarray*}
I_{\hat{B}}(g)=\inf\Big\{\frac12 \int_0^T ||\dot{f}(s)||^2 \, ds  : \hat f=g, \,\, f\in H_0^{1,p}[0,T]\Big\},
\end{eqnarray*}
with the understanding $ I_{\hat{B}}(g)= +\infty $ if the set
is empty. Thanks to  Proposition \ref{prop:LDPcc-uncorr}, $(X^{n,\varphi})_n$ where 
$$X^{n,\varphi}=X^n|\big(\hat{B}^n(t)=\varphi(t)\quad 0\leq t\leq T\big)$$ satisfies the LDP continuity condition. Therefore we have that the family $(\hat{B}^n,X^n)_{n\in\N}$ satisfies the hypotheses of Chaganty's Theorem and then  enjoys the WLDP with the rate function given by (\ref{eq:rate-function-couple}). Always by Chaganty's Theorem the family of processes $(X^n)_{n\in\N}$ satisfies a LDP  with the  speed $\ep_n^{-2}$ and the rate function
\begin{eqnarray*}I_X(x)
&=& \inf_{\varphi\in  \cl C_0 ^p}I(\varphi,x)
=\inf_{\varphi\in  \cl C_0 ^p}\big(I_{\hat{B}}(\varphi) +J(x|{\varphi})\big).\end{eqnarray*}
The first term in the sum is infinite if does not exist   $f\in H_0^{1,p}[0,T]$ such that $\varphi=\hat f$. Therefore we have 
that  $I_X(\,\cdot\,)$ is given by equation
(\ref{eq:log-price-rate-function-uncorr}).\cvd

We proved that the family of log-price processes $(X^n)_{n\in\mathbb{N}}$ satisfies a LDP with the  speed $\ep_n^{-2}$ and the rate  function $I_X(\,\cdot\,)$. The proof that  $I_X(\,\cdot\,)$ is  a good rate function follows from  Lemma 2.6 in \cite{Cha}.
\begin{proposition}\label{prop:I-X-good}
	The rate function $I_X(\,\cdot\,)$ is a good rate function.
\end{proposition}

\proof

It is enough to show (see also \cite{Cel-Pac} for the details), that
	$$\bigcup_{\varphi\in K_1}\{x\in { \cl C_0 }^d:J(x|\varphi)\leq L\}$$
	is a compact subset of ${ \cl C_0 }^d$ for any $L\geq 0$ and for any compact set $K_1\subset { \cl C_0 }^p$.
Let $K_1$ be a compact subset of ${ \cl C_0 }^p$; for $\varphi\in K_1$ define
	$$A_\varphi^L=\{x\in { \cl C_0 }^d:J(x|\varphi)\leq L\}=\{ x\in H_0^{1,d}[0,T]: J(x|\varphi)\leq L \}.$$
	$A_\varphi^L$ is a compact subset of ${ \cl C_0 }^d$.  We want to show that every sequence in $\displaystyle\bigcup_{\varphi\in K_1}A_\varphi^L$ has a convergent subsequence.
Let $\displaystyle(x_n)_{n\in\mathbb{N}}\subset\bigcup_{\varphi\in K_1}A_\varphi^L$, then, for every $n\in\mathbb{N}$, there exists $\varphi_n\in K_1$ such that $x_n\in A_{\varphi_n}^L$ (i.e. $J(x_n|\varphi_n)\leq L$). Since $(\varphi_n)_{n\in\mathbb{N}}\subset K_1$, up to a subsequence, we can suppose that $\varphi_n\to\varphi$ in ${ \cl C_0 }^p$, as $n\to+\infty$, with $\varphi\in K_1$.
Now, we claim that there exists a constant $N>0$ such that, for every $n\in\mathbb{N}$,
	$J(x_n|\varphi)\leq N.$
	Note that, 
	\begin{equation*}\begin{split}
		  J(x_n|\varphi)&=\Gamma\Big(x_n-\int_0^\cdot\mu(\varphi(t))\,dt\Big| a^{-1}(\varphi)\Big)\\&=\begin{cases} \frac12\int_0^T(\dot{x}_n(t)-\mu(\varphi(t)))^Ta^{-1}(\varphi(t))(\dot{x}_n(t)-\mu(\varphi(t)))\,dt&x\in H_0^{1,d}[0,T]\\
			+\infty&otherwise.
		\end{cases}\end{split}\end{equation*}
Therefore adding and subctracting  the term $\int_0^\cdot\mu(\varphi_n(t))dt$ and taking into account  Remark \ref{rem:Gamma}$(ii)$, for every $n\in\mathbb{N}$, we have,
$$\displaylines{
	J(x_n|\varphi)=\Gamma\Big(x_n-\int_0^\cdot\mu(\varphi(t))dt\Big| a^{-1}(\varphi)\Big)=\Gamma\Big(x_n-\int_0^\cdot\mu(\varphi_n(t))dt+\int_0^\cdot\mu(\varphi_n(t))-\int_0^\cdot\mu(\varphi(t))dt\Big| a^{-1}(\varphi)\Big)\cr
	\leq2\Gamma\Big(x_n-\int_0^\cdot\mu(\varphi_n(t))dt\Big| a^{-1}(\varphi)\Big)+2\Gamma\Big(\int_0^\cdot(\mu(\varphi_n(t))-\mu(\varphi(t)))dt\Big| a^{-1}(\varphi)\Big).}
	$$
	Thanks to Lemma \ref{lemma:tech3}, there exists a constant $M>1$ such that, for every $n\in\mathbb{N}$, the matrix $Ma^{-1}(\varphi_n)-a^{-1}(\varphi)$ is  positive definite; thanks to Remark \ref{rem:unif-def-pos},  there exists a constant $\beta_\varphi>0$ such that $\beta_\varphi I-a^{-1}(\varphi) $  is positive definite, therefore by using  Remark \ref{rem:Gamma}$(i)$ we have 
	$$\displaylines{
	J(x_n|\varphi)\leq 2MJ(x_n|\varphi_n)+\beta_\varphi\int_0^T(\mu(\varphi_n(t))-\mu(\varphi(t)))^T(\mu(\varphi_n(t))-\mu(\varphi(t)))dt\leq N,}$$
 since $J(x_n|\varphi_n)\leq L$ and the integral is bounded.
Then, since $A_\varphi^N=\{x\in { \cl C_0 }^d: J(x|\varphi)\leq N\}$ is a compact subset of ${ \cl C_0 }^d$, up to a subsequence, we can suppose that $x_n\to x\in A_\varphi^N$ in ${ \cl C_0 }^d$, as $n\to+\infty$.
Moreover, from the semicontinuity of $J(\cdot|\cdot)$, we have
	$$J(x|\varphi)\leq \liminf_{n\to +\infty}J(x_n|\varphi_n)\leq L,$$
	that implies $x\in A_\varphi^L$.\cvd

\section{The correlated model} \label{sect:correlated}

The dynamic of the asset price process $(S(t))_{t\in[0,T]}=(S_1(t),\ldots,S_d(t))_{t\in[0,T]}$ is modeled by \eqref{eq:general-model}.
 
Let $Z_i(t)=\log S_i(t)$ be the $i$-th log-price process, then we have
\begin{equation}\label{eq:log-price-corr}\begin{array}{c}
	Z_i(t)=
X_i(t) -\frac12\int_0^t\sum_{\ell=1}^p\tilde{\sigma}_{i\ell}(\hat{B}(s))^2 ds+\int_0^t\sum_{\ell=1}^p\tilde{\sigma}_{i\ell}(\hat{B}(s))\,dB_\ell(s)\phantom{asdadrrwrwr}
\end{array}\end{equation}
for every $1\leq i\leq d$ and $0\leq t \leq T$, where,   $(X(t))_{t\in[0,T]}$ is the  log-price process defined in the previous section.
Now,  consider the following scaled version of the stochastic differential equation in (\ref{eq:general-model})
$$
	\begin{cases}
	 \frac{dS_i^n(t)}{S_i^n(t)}=\mu_i(\hat{B}^n(t))\,dt+\ep_n\Big(\sum_{\ell=1}^p\tilde{\sigma}_{i\ell}(\hat{B}^n(t))\,dB_\ell(t)+\sum_{j=1}^d\sigma_{ij}(\hat{B}^n(t))\,dW_j(t)\Big)\\
		S_i^n(0)=1
	\end{cases}
$$
for every $1\leq i\leq d$ and $0\leq t\leq T$,  where $(\hat{B}^n)_{n\in\N}$
is a family of multidimensional Volterra processes satisfying  Assumption \ref{ass:LDP-B-hatB}.
 Then, the $i$-th log-price process in the scaled model is
\begin{equation}\label{eq:scaled-log-price-corr}
	Z_i^n(t)=X_i^n(t)-\frac12\ep_n^2\int_0^t\sum_{\ell=1}^p\tilde{\sigma}_{i\ell}(\hat{B}^n(s))^2\,ds+\ep_n\int_0^t\sum_{\ell=1}^p\tilde{\sigma}_{i\ell}(\hat{B}^n(s))\,dB_\ell(s)
\end{equation}
where, for every $n\in\N$, $(X^n(t))_{t\in[0,T]}$ is the scaled log-price process defined in the previous section.

In what follows, we will want to prove a sample path LDP for the family of processes  $(Z^n)_{n\in\N}$, defined in (\ref{eq:scaled-log-price-corr}). The study of the correlated model is more complicated than the previous one.
In fact, in this case, we should also study the behavior of the family of processes $$\Big(\Big(V_i^n(t)-\frac12\ep_n^2\int_0^t\sum_{\ell=1}^p\tilde{\sigma}_{i\ell}(\hat{B}^n(s))^2\,ds\Big)_{t\in[0,T]}\Big)_{n\in\N},$$
where
\begin{equation}\label{eq:process-V}
	V_i^n(t)=\ep_n\int_0^t\sum_{\ell=1}^p\tilde{\sigma}_{i\ell}(\hat{B}^n(s))\,dB_\ell(s)
\end{equation}
for every $1\leq i\leq d$, $0\leq t\leq T$.
Notice that this process depends on the couple $(\ep_n B, \hat{B}^n)$, but we can't directly apply  Chaganty's Theorem to the family
$$ ((\ep_n B, \hat{B}^n), Z^{n})_{n \in \mathbb{N}} $$
since
$ V^n $ cannot be written as a continuous function of $ (\ep_n B, \hat{B}^n) $ and so the LDP continuity condition is not fulfilled. To overcome this problem, as in \cite{Cel-Pac}, we  introduce a new family of processes $ (Z^{n,m})_{n \in \mathbb{N}} $, where for every $ m \geq 1 $, $V^n$
is replaced by a suitable continuous function  of $ (\ep_n B, \hat{B}^n) $.
Thanks to the results obtained in the previous section, we  prove that the hypotheses of Chaganty's Theorem are fulfilled for the family $ ((\ep_n B, \hat{B}^n), Z^{n,m})_{n \in \mathbb{N}} $. Then, for every $ m \geq 1 $, $ (Z^{n,m})_{n \in \mathbb{N}} $ satisfies a LDP with a certain good rate function $ I^m$ (Section 4.1). Then, proving that the family $ ((Z^{n,m})_{n \in \mathbb{N}})_{m\geq 1} $ is an exponentially good approximation (see Definition \ref{def:exp-approx}) of $ (Z^{n})_{n \in \mathbb{N}} $, we  obtain a LDP for the family $ (Z^{n})_{n \in \mathbb{N}} $ with the good rate function obtained in terms of the $ I^m $'s (Section 4.2). Finally (in Section 4.3) we give an explicit expression for the rate function (not in terms of the $ I^m $'s).

Suppose $\sigma\in {\cl C}^{d\times d}$ and $\mu\in {\cl C}^{d}$ satisfy Assumption \ref{ass:hp-sigma-mu-I}. In this section, we need some more hypotheses on the coefficients.


\begin{assumption}\label{ass:hp-sigma-tilde}
	 $\tilde{\sigma}_{i\ell}:\R^p\to\R$ are locally $\omega$-continuous functions (see Definition 6.1 in \cite{Cel-Pac}) for $1\leq i\leq d$ and $1\leq \ell\leq p$.
\end{assumption}
\begin{assumption}\label{ass:hp-mu-sigma-II}
\rm
		 There exist constants $\alpha, M_1,M_2>0, $ such that for every $i,j=1,\ldots,d$, $1\leq \ell\leq p$,
		$$ |\tilde\sigma_{i\ell}(x)|+|\sigma_{ij}(x)|+ |\mu_i(x)|\leq M_1+ M_2 \,||x||^\alpha, \quad  x\in \R^p.$$
\end{assumption}

\begin{remark}\label{rem:eigenvalue-estimate}\rm
Under  Assumption \ref{ass:hp-mu-sigma-II}
if $\lambda_{max}(x)$ is the maximum eigenvalue of the matrix $a(x)$, then there exists a constant $M>0$ such that
 $\lambda_{max}(x)\leq M ||x||^{ 2\alpha}$ and therefore
$\frac 1{\lambda_{max}(x)}\geq \frac 1{M ||x||^{2 \alpha}}$. Notice that $\frac 1{\lambda_{max}(x)}$ is the minimum eigenvalue of the matrix
$a^{-1}(x)$.

\end{remark}
The main result of this section is the following theorem. Also in this section we procede as in \cite{Cel-Pac}.
\begin{theorem}\label{th:LDP-main}
	Suppose that  Assumptions \ref{ass:LDP-B-hatB}, \ref{ass:hp-sigma-mu-I},  \ref{ass:hp-sigma-tilde} and \ref{ass:hp-mu-sigma-II} are fulfilled. Then, the family of processes $(Z^n)_{n\in\mathbb{N}}$ satisfies a LDP with the  speed $\ep_n^{-2}$ and the good rate function
	\begin{equation}\label{eq:rate-main}
		\mathcal{I}_Z(x)=\begin{cases}\displaystyle\inf_{f\in H_0^{1,p}[0,T]}\mathcal{I}((f,\hat{f}),x)\quad &x\in H_0^{1,d}[0,T]\\
			+\infty&otherwise
		\end{cases}
	\end{equation}
	where for every $f\in H_0^{1,p}[0,T]$,
	\begin{equation*}
		\mathcal{I}((f,\hat{f}),x)= \frac12\|f\|_{H_0^{1,p}[0,T]}^2
		+J(x-\Phi(f,\hat{f})|\hat{f})
	\end{equation*}
	where $J(\,\cdot\,|\,\cdot\,)$ is defined in \eqref{eq:cond-rate-function-uncorr} and  $\Phi$ is defined in \eqref{eq:Phi}.
\end{theorem}

\subsection{LDP for the Approximating Families}

In this section we suppose that Assumptions \ref{ass:hp-sigma-mu-I} and \ref{ass:hp-sigma-tilde} are fulfilled.
 First we  define the analogous of the function $\Psi^m$  in equation (22)  in \cite{Cel-Pac}.

For every $m\geq 1$, let us define the function $\Phi^m: \cl C_0 ^p\times  \cl C_0 ^p\to {\cl C}_0^d$ as follows: for $(f,g)\in  \cl C_0 ^p\times  \cl C_0 ^p$, $t\in[0,T]$ and $i=1,\ldots,d$,
\begin{equation}\label{eq:Phi^m}
\Phi_i^m(f,g)(t)=\sum_{\ell=1}^p\Psi_{i\ell}^m(f,g)(t),\end{equation}
where
$$
	\Psi_{i\ell}^m(f,g)(t)=\sum_{k=0}^{\bigl\lfloor\frac{mt}{T}\bigr\rfloor -1}\tilde{\sigma}_{i\ell}\Big(g\Big(\frac{k}{m}T\Big)\Big)\Big[f_\ell\Big(\frac{k+1}{m}T\Big)-f_\ell\Big(\frac km T\Big) \Big]+ \tilde{\sigma}_{i\ell}\Big(g\Big(\Bigl\lfloor\frac{mt}{T}\Bigr\rfloor\frac Tm\Big)\Big)\Big[f_\ell(t)-f_\ell\Big(\Bigl\lfloor\frac{mt}{T}\Bigr\rfloor\frac Tm\Big)\Big]$$
and $\lfloor\cdot\rfloor$ is the usual floor function.
Let us collect some properties of the function $\Phi^m(\cdot,\cdot)$ in the following remark.
\begin{remark}\label{rem:Phi^m}\rm
 It is clear that $\Phi^m(\cdot,\cdot)$ is a continuous function on $ \cl C_0 ^p\times  \cl C_0 ^p$ (where we are using the sup norm topology for both arguments).

Furthermore
for every $(f,g)\in H_0^{1,p}[0,T]\times  \cl C_0 ^p$ the function $\Phi^m(\cdot,\cdot)$ can be written in the following way:
	$$
		\Phi_i^m(f,g)(t)=\sum_{\ell=1}^p\int_0^t\tilde{\sigma}_{i\ell}\Big(g\Big( \Bigl\lfloor\frac{ms}{T}\Bigr\rfloor\frac Tm\Big)\Big)\dot{f}_\ell(s)\,ds,
	$$
for every $1\leq i\leq d$ and $0\leq t\leq T$. Then
 thanks to the Cauchy-Schwarz inequality we have that
	for every $f\in H_0^{1,p}[0,T]$, there exists $M>0$, depending on $f$, such that,
	$$\int_0^T\|\dot{\Phi}^m(f,\hat{f})(t)\|^2dt\leq M\|f\|_{H_0^{1,p}[0,T]}^2.$$
\end{remark}

\medskip

For every $m\geq 1$, let us define the new family of processes $((Z^{n,m}(t))_{t\in[0,T]})_{n\in\N}$, where
\begin{equation*}
\begin{split}
	Z_i^{n,m}(t)&=\int_0^t\Big(\mu_i(\hat{B}^n(s))-\frac12\ep_n^2\sum_{j=1}^d\sigma_{ij}(\hat{B}^n(s))^2-\frac12\ep_n^2\sum_{\ell=1}^p\tilde{\sigma}_{i\ell}(\hat{B}^n(s))^2\Big) ds\phantom{dsfsdsfsdfsdfsdfds}\\
	&\phantom {dzxczxczxfsfsdfsdfadfsdfs}+\Phi_i^m(\ep_nB,\hat{B}^n)(t)+\ep_n\int_0^t\sum_{j=1}^d\sigma_{ij}(\hat{B}^n(s))\,dW_j(s),
\end{split}
\end{equation*}
where $\Phi_i^m(\cdot,\cdot)$ is defined in (\ref{eq:Phi^m}), for every $1\leq i\leq d$ and $0\leq t\leq T$.
For every $m\geq 1$, we will prove a LDP  for the family of processes $(Z^{n,m})_{n\in\N}$. 
More precisely we prove the following theorem.
\begin{theorem}\label{th:LDP-approx}
	Suppose that  Assumptions \ref{ass:hp-sigma-mu-I} and  \ref{ass:hp-sigma-tilde} are fulfilled. For every $m\geq 1$, a LDP with the speed $\ep_n^{-2}$ and the good rate function $I_Z^m(\,\cdot\,)$ given by (\ref{eq:rate-function-I-Z^m}) holds  for the family  $(Z^{n,m})_{n\in\N}$.
\end{theorem}For this purpose we will check that hypotheses of Chaganty's Theorem hold for the family of processes
$((\ep_nB,\hat{B}^n),Z^{n,m})_{n\in\N}.$
By Assumption \ref{ass:LDP-B-hatB}, we already know that $((\ep_nB,\hat{B}^n))_{n\in\N}$ satisfies a LDP  with the  speed $\ep_n^{-2}$ and the good rate function $I_{B,\hat{B}}(\cdot,\cdot)$ given by (\ref{eq:rate-function-B-hatB}). Fix $m\geq 1$ and $(f,g)\in  \cl C_0 ^p\times  \cl C_0 ^p$. Our next goal is to prove that the family of conditional processes
$$Z^{n,m,(f,g)}=Z^{n,m}|\big(\ep_nB(t)=f(t), \hat{B}^n(t)=g(t)\quad 0\leq t\leq T\big)$$
satisfies the LDP continuity condition (condition \textit{(ii)} of Chaganty's Theorem.

For every $(f,g)\in  \cl C_0 ^p\times  \cl C_0 ^p$ and $t\in[0,T]$, in law, for every $1\leq i\leq d$, we have
\begin{equation}\label{eq:Z^m-scaled-cond}Z_i^{n,m,(f,g)}(t)
	=\bar{X}_i^{n,g}(t)+\Phi_i^m(f,g)(t)\end{equation}
  with $\bar{X}^{n,g}$ defined as
\begin{eqnarray}\label{eq:barX-cond}
	\bar{X}_i^{n,g}(t)=X_i^{n,g}(t)-\frac12\ep_n^2\int_0^t\sum_{\ell=1}^p\tilde{\sigma}_{i\ell}(g(s))^2\,ds
\end{eqnarray}
and  $X^{n,g}$ is the conditional log price process in the uncorrelated model defined in (\ref{eq:logprice-uncorr-scaled-cond}).

The following proposition proves the  LDP continuity condition for the family  $(Z^{n,m,(f,g)})_{n\in\N}$.
\begin{proposition} \label{prop:LDPcc-corr}
Fix  $(f,g)\in  \cl C_0 ^p\times  \cl C_0 ^p$. Then, for every $m\geq 1$, the sequence of the conditional processes  $(Z^{n,m,(f,g)})_{n\in\N}$  satisfies the LDP continuity condition  with the  speed $ \ep_n^{-2}$ and the rate function
		  \begin{equation}\label{eq:rate-function-J^m}
		\mathcal{J}^m(x|(f,g))=J(x-\Phi^m(f,g)|g)
	\end{equation}
where $J(\,\cdot\,|g)$ is given by (\ref{eq:cond-rate-function-uncorr}). Notice that $\mathcal{J}^m(x|(f,g))$ is finite if and only if $x-\Phi^m(f,g)\in H_0^{1,d}[0,T]$.
 \end{proposition}
\proof
	
$(a)$ 	
Fix $(f,g)\in  \cl C_0 ^p\times  \cl C_0 ^p$. We have to prove that the sequence of the conditional processes  $(Z^{n,m,(f,g)})_{n\in\N}$  satisfies the LDP continuity condition  with the  speed $ \ep_n^{-2}$ and the rate function defined in (\ref{eq:rate-function-J^m}).
It is not hard to prove that 	for every $g\in  \cl C_0 ^p$ the family of  processes $(\bar{X}^{n,g})_{n\in\N}$ satisfies the same LDP as  $({X}^{n,g})_{n\in\N}$ then
thanks to  equation (\ref{eq:Z^m-scaled-cond}) this  is a simple applications of the contraction principle.

\medskip

$(b)$ 	
	Let $((f_n,g_n))_{n\in\N}\subset  \cl C_0 ^p\times  \cl C_0 ^p$ be such that
	$(f_n,g_n)\to(f,g)$
	in $ \cl C_0 ^p\times  \cl C_0 ^p$, as $n\to+\infty$.
We have to prove that  the family of processes $(Z^{n,m,(f_n,g_n)})_{n\in\N}$, for  every $m\geq 1$,
	 satisfies a LDP   with the  speed $\ep_n^{-2}$ and the good rate function $\mathcal{J}^m(\,\cdot\,|(f,g))$ defined in (\ref{eq:rate-function-J^m}).

If $g_n\to g$ in $ \cl C_0 ^p$, as $n\to+\infty$, immediately follows that
the family of processes $((\bar{X}^{n,g_n})_{n\in\N}$ satisfies a LDP  with the  speed $\ep_n^{-2}$ and the good rate function $J(\,\cdot\,|g)$.
	Combining this with contraction principle, we have that, for every $m\geq 1$, the family of processes
	$$(\bar{X}^{n,g_n}+\Phi^m(f,g))_{n\in\N}$$
	satisfies a LDP with the  speed $\ep_n^{-2}$ and the good rate function $\mathcal{J}^m(\,\cdot\,|(f,g))$, defined in (\ref{eq:rate-function-J^m}). Furthermore, for every $m\geq 1$,
	$\Phi^m(f_n,g_n)\to\Phi^m(f,g)$
	in $ (\cl C_0 ^p)^2$, as $n\to+\infty$, since $\Phi^m$ is a continuous function.
Therefore, the families $(\bar{X}^{n,g_n}+\Phi^m(f_n,g_n))_{n\in\N}$ and $(\bar{X}^{n,g_n}+\Phi^m(f,g))_{n\in\N}$ are exponentially equivalent
(see Defintion 4.2.10 in \cite{Dem-Zei}) and the statement is proved.

\medskip

$(c)$
We have to prove the lower semi-continuity of $\mathcal{J}^m(\,\cdot\,|(\cdot,\cdot))$ as a function of $(x,(f,g))\in \cl C_0 ^d$ $\times \big( \cl C_0 ^p\times  \cl C_0 ^p\big)$. That is
	if the sequence $((f_n,g_n),x_n)_{n\in\N}$ is such that
	$$((f_n,g_n),x_n)\to((f,g),x)$$
	in $ \cl C_0 ^p$  $\times  \cl C_0 ^p\times \cl C_0 ^d$, as $n\to+\infty$, then, for every $m\geq 1$,
	$$\liminf_{n\to +\infty}\mathcal{J}^m(x_n|(f_n,g_n))\geq \mathcal{J}^m(x|(f,g)).$$

	For every $m\geq 1$, $\Phi^m(\cdot,\cdot)$ is continuous on $ \cl C_0 ^p\times  \cl C_0 ^p$, therefore if $((f_n,g_n),x_n)\to((f,g),x)$ in $\big( \cl C_0 ^p\times  \cl C_0 ^p\big)\times \cl C_0 ^d$, as $n\to+\infty$, then
	$$x_n-\Phi^m(f_n,g_n)\to x-\Phi^m(f,g)$$
	in $\cl C_0 ^d$, as $n\to+\infty$. Then, by the lower semi-continuity of $J(\,\cdot\,|\,\cdot\,)$,
	\begin{eqnarray*}
		\liminf_{n\to +\infty}\mathcal{J}^m(x_n|(f_n,g_n))=\liminf_{n\to +\infty}J(x_n-\Phi^m(f_n,g_n)|g_n)
		\geq J(x-\Phi^m(f,g)|g)=\mathcal{J}^m(x|(f,g))
	\end{eqnarray*}
	and this completes the proof that, for every $(f,g)\in  \cl C_0 ^p\times  \cl C_0 ^p$ the family of processes $(Z^{n,m,(f,g)})_{n\in\N}$ satisfies the LDP continuity condition for every $m\geq 1$. \cvd

\begin{proposition}
	Suppose $\sigma$, $\mu$ and $\tilde{\sigma}$ satisfy Assumptions \ref{ass:hp-sigma-mu-I} and \ref{ass:hp-sigma-tilde}.
	Then, for every $m\geq 1$, the family $((\ep_nB,\hat{B}^n),Z^{n,m})_{n\in\N}$ satisfies a WLDP with the  speed $\ep_n^{-2}$ and the rate function
	$$
		I^m((f,g),x)=I_{B,\hat{B}}(f,g)+J(x-\Phi^m(f,g)|g)
$$
	for $x\in { \cl C_0 }^d$ and $(f,g)\in  \cl C_0 ^p\times  \cl C_0 ^p$, where $I_{B,\hat{B}}(\cdot,\cdot)$ and $J(\,\cdot\,|\,\cdot\,)$ are defined in (\ref{eq:rate-function-B-hatB}) and (\ref{eq:cond-rate-function-uncorr}), respectively.
Moreover, the family of processes $(Z^{n,m})_{n\in\N}$ satisfies a LDP  with the  speed $\ep_n^{-2}$ and the rate function
	\begin{equation}\label{eq:rate-function-I-Z^m}
	  I_Z^m(x)=\begin{cases}
		\displaystyle	\inf_{f\in H_0^{1,p}[0,T]}\Bigl\{ \frac12\|f\|_{H_0^{1,p}[0,T]}^2+J(x-\Phi^m(f,\hat{f})|\hat{f})\Bigr\}\quad&x\in H_0^{1,d}[0,T]\\
			\phantom{dd}+\infty,&otherwise
		\end{cases}
	\end{equation}
where $\hat{f}$ is defined in (\ref{eq:hatf}).
	
\end{proposition}

\proof
	Thanks to Assumption \ref{ass:LDP-B-hatB} and  Proposition \ref{prop:LDPcc-corr}, the family $((\ep_nB,\hat{B}^n),Z^{n,m})_{n\in\N}$ satisfies the hypothesis of Chaganty's Theorem. Therefore $(Z^{n,m})_{n\in\N}$, satisfies a LDP   with the  speed $\ep_n^{-}2$ and the rate function
	$$I_Z^m(x)=\inf_{\{(f,g)\in  \cl C_0 ^p\times  \cl C_0 ^p\}}\Bigl\{ I_{B,\hat{B}}(f,g)+\mathcal{J}^m(x|(f,g)) \Bigr\}$$
	for $x\in \cl C_0 ^d$. Combining the expressions of $I_{B,\hat{B}}(\cdot,\cdot)$ and $\mathcal{J}^m(\cdot|(\cdot,\cdot))$, the rate function $I_Z^m(\,\cdot\,)$ is given by (\ref{eq:rate-function-I-Z^m}).
\cvd

We have proved that, for every $m\geq 1$, the family of processes $(Z^{n,m})_{n\in\N}$ satisfies a LDP with the  speed $\ep_n^{-2}$ and the rate  function $I_Z^m(\,\cdot\,)$. Now, we  want to prove that, for every $m\geq 1$, the rate function $I_Z^m(\,\cdot\,)$ is actually a good rate function.
The proof of the following proposition is the same as Proposition 6.11 in \cite{Cel-Pac}. We give the details since some estimates are not immediate.
\begin{proposition}\label{prop:I-Z^m-good}
	For every $m\geq 1$, the rate function $I_Z^m(\,\cdot\,)$ is a good rate function.
\end{proposition}
\proof
As in Proposition \ref{prop:I-X-good} it is enough to show that
	 the set
	$$\bigcup_{(f,g)\in K_1}\{x\in { \cl C_0 }^d:\mathcal{J}^m(x|(f,g))\leq L\}=\bigcup_{(f,g)\in K_1} A_{(f,g)}^L,$$
	is a compact subset of ${ \cl C_0 }^d$ for any $L\geq 0$ and for any level set $K_1$ of the rate function $I_{B,\hat{B}}(\cdot,\cdot)$.
Let $K_1$ be a level set of $I_{B,\hat{B}}(\cdot,\cdot)$ and  $L\geq 0$. If  $(f,g)\in K_1$,  then $g=\hat f$ and
	$$A_{(f,g)}^L=A_{(f,\hat{f})}^L=\{x\in H_0^{1,d}[0,T]:\mathcal{J}^m(x|(f,\hat{f}))\leq L\}.$$
	
	For every $(f,\hat{f})\in K_1$, $A_{(f,\hat{f})}^L$ is a compact set of ${ \cl C_0 }^d$, since $\mathcal{J}^m(\,\cdot\,|(f,\hat{f}))$ is a good rate function. We want to show that every sequence in $\displaystyle\bigcup_{(f,\hat{f})\in K_1}A_{(f,\hat{f})}^L$ has a convergent subsequence.
	 Let $\displaystyle(x_n)_{n\in\N}\subset\bigcup_{(f,\hat{f})\in K_1}A_{(f,\hat{f})}^L$, then, for every $n\in\N$, there exists $(f_n,\hat{f}_n)\in K_1$ such that $x_n\in A_{(f_n,\hat{f}_n)}^L$ (that is $\mathcal{J}^m(x_n|(f_n,\hat{f}_n))\leq L$).  Since $((f_n,\hat{f}_n))_{n\in\N}\subset K_1$, up to a subsequence, we can suppose $(f_n,\hat{f}_n)\to(f,\hat{f})$ in $ \cl C_0 ^p\times  \cl C_0 ^p$, as $n\to+\infty$, with $(f,\hat{f})\in K_1$.
	 Now, we claim that there exists a constant $N>0$, such that, for every $n\in\N$,
	$$\mathcal{J}^m(x_n|(f,\hat{f}))=\Gamma\Big(x_n-\Phi^m(f,\hat{f})-\int_0^\cdot\mu(\hat{f}(t))dt\Big| a^{-1}(\hat{f})\Big)\leq  N.$$
	 From  \textit{(iii)} in Remark \ref{rem:Gamma} and Lemma \ref{lemma:tech3}, adding and subtracting $\displaystyle\Phi^m(f_n,\hat{f_n})+\int_0^\cdot\mu(\hat{f}_n(t))dt$, we have that, for every $n\in\N$,
	$$ \displaylines{
		\mathcal{J}^m(x_n|(f,\hat{f}))=\Gamma\Big(x_n-\Phi^m(f,\hat{f})-\int_0^\cdot\mu(\hat{f}(t))dt\Big| a^{-1}(\hat{f})\Big)
		\cr \leq 3M  \mathcal{J}^m(x_n|(f_n,\hat{f}_n))
		+3\Gamma\Big(\Phi^m(f_n,\hat{f}_n)-\Phi^m(f,\hat{f})\Big|a^{-1}(\hat{f})\Big)
		+3\Gamma\Big(\int_0^\cdot(\mu(\hat{f}_n(t))-\mu(\hat{f}(t)))dt\Big|a^{-1}(\hat{f})\Big)\cr
		\phantom{m}\leq3ML+3\Gamma\Big(\Phi^m(f_n,\hat{f}_n)-\Phi^m(f,\hat{f})\Big|a^{-1}(\hat{f})\Big) +3\Gamma\Big(\int_0^\cdot(\mu(\hat{f}_n(t))-\mu(\hat{f}(t)))dt\Big|a^{-1}(\hat{f})\Big).}$$

	The term $\displaystyle \Gamma\Big(\Phi^m(f_n,\hat{f}_n)-\Phi^m(f,\hat{f})\Big| a^{-1}(\hat{f})\Big)$ is bounded by a constant independent of $n$. Indeed, combining Remarks \ref{rem:unif-def-pos} and \ref{rem:Phi^m}, there exists a constant $M>0$ (it may also change from line to line) such that
	$$
		\Gamma\Big(\Phi^m(f_n,\hat{f}_n)-\Phi^m(f,\hat{f})\Big| a^{-1}(\hat{f})\Big)\leq M\Big(\|f_n\|_{H_0^{1,p}[0,T]}^2+\|f\|_{H_0^{1,p}[0,T]}^2\Big)
	$$
and this quantity is  bounded, since $(f,\hat{f}), (f_n,\hat{f}_n)\in K_1$, that is a level set of $I_{B,\hat{B}}(\cdot,\cdot)$.

Furthermore, reasoning as in Proposition \ref{prop:I-X-good}, we have
	$$\Gamma\Big(\int_0^\cdot(\mu(\hat{f}_n(t))-\mu(\hat{f}(t)))dt\Big|a^{-1}(\hat{f})\Big)\leq M.$$
	 Therefore there exists a constant $N>0$ such that $\mathcal{J}^m(x_n|(f,\hat{f}))\leq N$, that implies $x_n\in A_{(f,\hat{f})}^N$, for every $n\in\N$.  Then, since $A_{(f,\hat{f})}^N=\{x\in \cl C_0 ^d: \mathcal{J}^m(x|(f,\hat{f}))\leq N\}$ is a compact subset of $\cl C_0 ^d$, up to a subsequence, we can suppose that $x_n\to x \in A_{(f,\hat{f})}^N$ in ${ \cl C_0 }^d$, as $n\to+\infty$.  Moreover, from the lower semicontinuity of $\mathcal{J}^m(\cdot|(\cdot,\cdot))$,
	$$\mathcal{J}^m(x|(f,\hat{f}))\leq \liminf_{n\to +\infty}\mathcal{J}^m(x_n|(f_n,\hat{f}_n))\leq L,$$
	that implies $x\in A_{(f,\hat{f})}^L$.
\cvd

\subsection{LDP for the log-price process}
Theorem \ref{th:LDP-approx} provides a LDP for the families $ (Z^{n,m})_{n \in \mathbb{N}} $ for every $ m\geq1, $ but our goal is to get a LDP for the family $ (Z^{n})_{n \in \mathbb{N}} $.

We have to  prove that the sequence of processes $ ((Z^{n,m})_{n \in \mathbb{N}})_{m \geq 1} $ is an exponentially good approximation of $ (Z^{n})_{n \in \mathbb{N}} $.

Our next goal is to get a LDP for the family of log-prices $(Z^n)_{n\in\mathbb{N}}$, defined in (\ref{eq:scaled-log-price-corr}). We will prove that, for every $m\geq 1$, the families of processes $(Z^{n,m})_{n\in\mathbb{N}}$  are exponentially good approximations of $(Z^n)_{n\in\mathbb{N}}$.
\begin{proposition}\label{prop:good-approx}
	The families  $(Z^{n,m})_{n\in\mathbb{N}}$ are exponentially good approximations of $(Z^n)_{n\in\mathbb{N}}$, for every $m\geq 1$.
\end{proposition}
\proof
	We need to check that, for every $\delta>0$,
	$$\lim_{m\to +\infty}\limsup_{n\to+\infty}\ep_n^2\log\mathbb{P}\big(\| Z^{n,m}-Z^n\|_\infty>\delta\big)=-\infty.$$
	Now,
	$$\|Z^{n,m}-Z^n\|_\infty=\|V^n-\Phi^m(\ep_nB,\hat{B}^n)\|_\infty=\sup_{t \in [0,T]}\|V^n(t)-\Phi^m(\ep_nB,\hat{B}^n)(t)\|$$
	with $V^n$ and $\Phi^m(\cdot,\cdot)$ defined in (\ref{eq:process-V}) and (\ref{eq:Phi^m}), respectively.
Therefore	$$
		\big\{\| Z^{n,m}-Z^n\|_\infty>\delta\big\}\subset
		\bigcup_{i=1}^d\biggl\{ \sup_{t \in [0,T]}|V_i^n(t)-\Phi_i^m(\ep_nB,\hat{B}^n)(t)|>\delta/\sqrt{d} \biggr\}
	$$
	that implies
	\begin{eqnarray*}
		&&\lim_{m\to +\infty}\limsup_{n\to+\infty}\ep_n^2\log\mathbb{P}\big(\| Z^{n,m}-Z^n\|_\infty>\delta\big)\\
		&&\leq\lim_{m\to +\infty}\limsup_{n\to+\infty}\ep_n^2\sum_{i=1}^d\log\mathbb{P}\bigg(\sup_{t\in[0,T]}|V_i^n(t)-\Phi_i^m(\ep_nB,\hat{B}^n)(t)|>\delta/\sqrt{d} \bigg)=-\infty
	\end{eqnarray*}
	where the last equality follows from Proposition 6.18 in \cite{Cel-Pac}.
\cvd

We have proved the following theorem.
\begin{theorem}\label{th:LDP-log-price-corr}
	Suppose that  Assumptions \ref{ass:LDP-B-hatB}, \ref{ass:hp-sigma-mu-I} and \ref{ass:hp-sigma-tilde} are fulfilled. Then, the family of processes $(Z^n)_{n\in\mathbb{N}}$, satisfies a LDP  with the  speed $\ep_n^{-2}$ and the good rate function $I_Z(\,\cdot\,)$ given by
	$${I}_Z(x)=\lim_{\delta\to0}\liminf_{m\to +\infty}\inf_{y\in B_\delta(x)}I_Z^m(y)=\lim_{\delta\to0}\limsup_{m\to +\infty}\inf_{y\in B_\delta(x)}I_Z^m(y),$$
	for $x\in {\cal C}_0^d$, with $B_\delta(x)=\{y\in {\cal C}_0^d:\|x-y\|_\infty<\delta\}$ and $I_Z^m(\,\cdot\,)$ defined in (\ref{eq:rate-function-I-Z^m}).
\end{theorem}

\subsection{Identification of the Rate Function}
In this section, we suppose that Assumptions \ref{ass:hp-mu-sigma-II}, \ref{ass:hp-sigma-tilde} and \ref{ass:hp-sigma-mu-I} are fulfilled.

Theorem \ref{th:LDP-log-price-corr}  provides a LDP for the family of processes $(Z^n)_{n\in\mathbb{N}}$ with a good rate function $I_Z(\,\cdot\,)$ obtained in terms of the $I_Z^m$s. Our next goal is to give an explicit expression to the rate function  $I_Z(\,\cdot\,)$.
Let us  define a function $\Phi: \cl C_0 ^p\times  \cl C_0 ^p\to { \cl C_0 }^d$ (the analogous of the function $\Psi$  in equation (35)  in \cite{Cel-Pac})
as
\begin{equation}\label{eq:Phi}
	\Phi_i(f,g)(t)=\begin{cases}
		\sum_{\ell=1}^p
	\Psi_{i\ell}(f,g)(t)&(f,g)\in\mathscr{H}_{(B,\hat{B})}^p\\
		0&otherwise
	\end{cases}
\end{equation}
where for every $1\leq i\leq d$, $1\leq \ell\leq p$
$$
	\Psi_{i\ell}(f,g)(t)=\begin{cases}\displaystyle
		\int_0^t\tilde{\sigma}_{i\ell}(\hat{f}(s))\dot{f}_\ell(s)\,ds\quad&(f,g)\in\mathscr{H}_{(B,\hat{B})}^p\\
		0&otherwise,
	\end{cases}
$$
and $\mathscr{H}_{(B,\hat{B})}^p$ is defined in (\ref{eq:RKHS-B-hatB}).
Since for $f_\ell\in H_0^1[0,T]$, we have that $\hat{f}_\ell\in { \cl C_0 }$, for every $1\leq \ell\leq p$, the function $\Phi(\cdot,\cdot)$ is finite on $  \cl C_0 ^p\times  \cl C_0 ^p$ and $ \Phi(f,\hat{f})$ is differentiable with a square integrable gradient, i.e. $ \Phi(f,\hat{f})\in H_0^{1,d}[0,T]$. Proceeding as in Remark \ref{rem:Phi^m},
from (\ref{eq:Phi}), we obtain that, for every $f\in H_0^{1,p}[0,T]$,
$$
\int_0^T\|\dot{\Phi}(f,\hat{f})(t)\|^2dt\leq M\|f\|_{H_0^{1,p}[0,T]}^2.
$$

Next lemma is the same as Lemma 2.13 in \cite{Gu2}. We give some details of the proof for the sake of completeness.
\begin{lemma}\label{lemma:Phi-Phi^m-convergence}
	For every $L>0$, let $D_L=\{f\in H_0^{1,p}[0,T] : \|f\|_{H_0^{1,p}[0,T]}^2\leq L\}$. Then one has,
	$$\displaystyle \lim_{m\to +\infty}\sup_{f\in D_L}\|\Phi(f,\hat{f})-\Phi^m(f,\hat{f})\|_{\infty}=0.$$
	
\end{lemma}
\proof As in  Lemma 22 in \cite{Gu1}, that still holds if the arguments are vectorial functions,  we have
\begin{equation}\label{eq:Gu1}
\lim_{m\to+\infty}\sup_{f \in D_{L}}\sup_{t \in [0,T]}\bigg|\tilde\sigma_{i\ell}(\hat{f}(t))-\tilde\sigma_{i\ell}\Big(\hat{f}\Big({\Big\lfloor \frac{mt}T\Big\rfloor}\frac Tm\Big)\Big)\bigg|=0,\end{equation}
for every $1\leq i\leq d$,	$1\leq \ell \leq p$.,
Furthermore,
for $1\leq i\leq d$,
$$\|\Phi_i(f,\hat{f})-\Phi_i^m(f,\hat{f})\|_{\infty}\leq\sum_{\ell=1}^p\|\Psi_{i\ell}(f,\hat{f})-\Psi_{i\ell}^m(f,\hat{f})\|_{\infty}.$$
Therefore from Lemma 6.22 in \cite{Cel-Pac} (the statement still holds when the argument are vectorial function), the  lemma is proved.
\cvd

We introduce the following functional
\begin{equation}\label{eq:I-Z-identification}
	\mathcal{I}_Z(x)=\begin{cases}\displaystyle\inf_{f\in H_0^{1,p}[0,T]}\mathcal{H}((f,\hat{f}),x)\quad &x\in H_0^{1,d}[0,T]\\
		+\infty&otherwise
	\end{cases}
\end{equation}
where for every $f\in H_0^{1,p}[0,T]$,
$$
	\mathcal{H}((f,\hat{f}),x)= \frac12\|f\|_{H_0^{1,p}[0,T]}^2+J(x-\Phi(f,\hat{f})|\hat{f})
$$
with $J(\,\cdot\,|\,\cdot\,)$ and $\Phi(\cdot,\cdot)$ defined in (\ref{eq:cond-rate-function-uncorr}) and  (\ref{eq:Phi}), respectively.
Now, we will enunciate some remarks and lemmas in order to prove that $I_Z(\cdot)={\cal I}_Z(\cdot)$.
\begin{remark}
	\rm For $ x\in H_0^{1,d}[0,T] $ we have,
		$$\mathcal{I}_Z(x)=\inf_{f \in H_0^{1,p}[0,T]}\mathcal{H}((f,\hat{f}),x)\leq \mathcal{H}((0,0),x)=\frac{1}2\int_0^T(\dot{x}(t)-\mu(0))^Ta^{-1}(0)(\dot{x}(t)-\mu(0))\,dt,$$
therefore
		$$\mathcal{I}_Z(x)= \inf_{f \in D_{C_x}}\mathcal{H}((f,\hat{f}),x)  $$
		where $C_x=\int_0^T(\dot{x}(t)-\mu(0))^Ta^{-1}(0)(\dot{x}(t)-\mu(0))\,dt$ and  $ D_{C_x}=\{f \in H_0^{1,p}[0,T]:\lVert f\rVert_{H_0^{1,p}[0,T]}^2\leq C_x\} $.
		Similarly, for $ x \in H_0^{1,d}[0,T], $
		for every $ m\geq 1,$ we have
 $$ I_{{Z}}^m(x)= \inf_{f \in D_{C_x}}\mathcal{H}_m((f,\hat{f}),x) $$
		 where, we recall,  $ I_{{Z}}^m(\cdot) $ is the rate function defined in (\ref{eq:rate-function-I-Z^m}) and
		$$\mathcal{H}_m((f,\hat{f}),x)=\frac12 \lVert f\rVert_{H_0^{1,p}[0,T]}^2+J(x-\Phi^m(f,\hat{f})|\hat{f}).$$
	\end{remark}

In order to prove that  $I_Z(\cdot) = \mathcal{I}_Z(\cdot) $, we have to verify  that the hypotheses of Proposition \ref{prop:identification} are fulfilled. We start by proving the convergence to $ \mathcal{I}_Z(\cdot) $ of the rate functions $ I_{{Z}}^m (\cdot)$'s.

\begin{lemma}\label{lemma:diff-Phi-punto-Phi^m-punto}
	For every  $x\in H_0^{1,d}[0,T]$, we have that
	$$\lim_{m\to +\infty}\sup_{f\in D_{C_x}}\|\dot{\Phi}(f,\hat{f})(t)-\dot{\Phi}^m(f,\hat{f})(t)\|_\infty^2= 0.$$
\end{lemma}
	\proof
	For $f\in D_{C_x}$, $M$ a positive constant, we have
	$$\displaylines{|\dot{\Phi}_i^m(f,\hat{f})(t)-\dot{\Phi}_i(f,\hat{f})(t)|^2
		\leq\Big(\sum_{\ell=1}^p|\dot{f}_\ell(t)|\sup_{t \in [0,T]}\Big|\tilde{\sigma}_{i\ell}\Big(\hat{f}\Big(\Bigl\lfloor\frac{mt}{T}\Bigr\rfloor\frac Tm \Big)\Big)-\tilde{\sigma}_{i\ell}(\hat{f}(t)) \Big|\Big)^2\cr
		\leq \Big(\sum_{\ell=1}^p|\dot{f}_\ell(t)|\Big)^2\Big(\sum_{\ell=1}^p\sup_{t \in [0,T]}\Big|\tilde{\sigma}_{i\ell}\Big(\hat{f}\Big(\Bigl\lfloor\frac{mt}{T}\Bigr\rfloor\frac Tm \Big)\Big)-\tilde{\sigma}_{i\ell}(\hat{f}(t)) \Big|\Big)^2\cr
		\leq M\sum_{\ell=1}^p|\dot{f}_\ell(t)|^2 \sum_{\ell=1}^p\sup_{t \in [0,T]}\Big|\tilde{\sigma}_{i\ell}\Big(\hat{f}\Big(\Bigl\lfloor\frac{mt}{T}\Bigr\rfloor\frac Tm \Big)\Big)-\tilde{\sigma}_{i\ell}(\hat{f}(t)) \Big|^2.}$$
	for every $1\leq i\leq d$. Therefore, the statements follows from equation (\ref{eq:Gu1}).
	\cvd

Now, we want to identify the rate function $I_Z(\,\cdot\,)$ for the family of processes $((Z^n(t))_{t\in[0,T]})_{n\in\mathbb{N}}$ with $\mathcal{I}_Z(\,\cdot\,)$. We start with the pointwise convergence (first point of the Proposition \ref{prop:identification}).
\begin{lemma}\label{lemma:pointwise-convergence}
	For every $x\in { \cl C_0 }^d$,
	$$\lim_{m\to +\infty}I_Z^m(x)=\mathcal{I}_Z(x),$$
	where $I_Z^m(\,\cdot\,)$ and $\mathcal{I}_Z(\,\cdot\,)$ are defined in (\ref{eq:rate-function-I-Z^m}) and (\ref{eq:I-Z-identification}), respectively.
\end{lemma}
\proof
	If $x\notin H_0^{1,d}[0,T]$, then, for every $m\geq 1$,  $I_Z^m(x)=\mathcal{I}_Z(x)=+\infty$. For $x\in H_0^{1,d}[0,T]$, we have
\begin{eqnarray*}
	|I_{{Z}}^m(x)- \mathcal{I}_Z(x)|&=&\bigg| \inf_{f \in D_{C_x}}\mathcal{H}_m((f,\hat{f}),x)-\inf_{f \in D_{C_x}}\mathcal{H}((f,\hat{f}),x) \bigg|
	\leq \sup_{f \in D_{C_x}}|\mathcal{H}_m((f,\hat{f}),x)-\mathcal{H}((f,\hat{f}),x)|.
	\end{eqnarray*}
	Straightforward computations show that
$$\displaylines{
	 \mathcal{H}_m((f,\hat{f}),x)-\mathcal{H}((f,\hat{f}),x)= \cr\Gamma(\Phi(f,\hat{f})-\Phi^m(f,\hat{f})|a^{-1}(\hat{f}))+
\int_0^T(\dot{\Phi}(f,\hat{f})(t)-\dot{\Phi}^m(f,\hat{f})(t))^Ta^{-1}(\hat{f}(t))(\dot{x}(t)-\mu(\hat{f}(t))-\dot{\Phi}(f,\hat{f})(t))\,dt.}$$
Now, thanks to Remark \ref{rem:unif-def-pos},
there exists a constant $M>0$ such that
	$$\displaylines{\sup_{f\in D_{C_x}}\Gamma(\Phi(f,\hat{f})-\Phi^m(f,\hat{f})|a^{-1}(\hat{f}))\leq M\sup_{f\in D_{C_x}} \int_0^T||\dot{\Phi}(f,\hat{f})(t)-\dot{\Phi}^m(f,\hat{f})(t) ||^2dt.}$$
The last term goes to zero thanks to Lemma \ref{lemma:diff-Phi-punto-Phi^m-punto}.
Moreover, 
for every $x\in H_0^{1,d}[0,T]$, thanks to Remark \ref{rem:hat-f-bounded} there exists a constant $A_x>0$ such that
	$$\sup_{f\in D_{C_x}}\|a^{-1}_{ij}(\hat{f})\|_{\infty}\leq A_x,\quad \sup_{f\in D_{C_x}}\|\mu_{i}(\hat{f})\|_{\infty}\leq A_x $$
	for every $1\leq i,j\leq d$.
Therefore thanks to the Cauchy-Schwarz inequality, we have
	$$  \displaylines{\sup_{f\in D_{C_x}}\Big|\int_0^T(\dot{\Phi}(f,\hat{f})(t)-\dot{\Phi}^m(f,\hat{f})(t))^Ta^{-1}(\hat{f}(t))(\dot{x}(t)-\mu(\hat{f}(t))-\dot{\Phi}(f,\hat{f})(t))\,dt\Big|\leq\cr
		A_x\sup_{f\in D_{C_x}}\sum_{i,j=1}^d\int_0^T|\dot{\Phi}_i(f,\hat{f})(t)-\dot{\Phi}_i^m(f,\hat{f})(t)||\dot{x}_j(t)-\mu_j(\hat{f}(t))-\dot{\Phi}_j(f,\hat{f})(t)|\,dt\leq\cr
	A_x\sup_{f\in D_{C_x}}\Big(\int_0^T\|\dot{\Phi}(f,\hat{f})(t)-\dot{\Phi}^m(f,\hat{f})(t)\|^2\,dt\Big)^{\frac12}\sup_{f\in D_{C_x}}\Big(\int_0^T\|\dot{x}(t)-\mu(\hat{f}(t))-\dot{\Phi}(f,\hat{f})(t)\|^2\,dt\Big)^{\frac12}.}$$	
It is not hard to prove that there exists $R_x>0$ such that,
$$\sup_{f\in D_{C_x}}\int_0^T\|\dot{x}(t)-\mu(\hat{f}(t))-\dot{\Phi}(f,\hat{f})(t)\|^2\,dt\leq R_x,$$
therefore, from Lemma \ref{lemma:diff-Phi-punto-Phi^m-punto},  one has
$\lim_{m\to +\infty}|I_Z^m(x)-\mathcal{I}_Z(x)|= 0$.
\cvd

It remains to show that  $x_m\overset{}{\underset{m \to +\infty}{\longrightarrow}}x  $ implies $  \liminf_{m\to+\infty} I_Z^m(x_m)\geq \mathcal{I}_Z(x)$
(second point of the Proposition \ref{prop:identification}). For this purpose we  need to prove that $ \mathcal{I}_Z(\cdot) $ is lower semicontinuous.

\begin{lemma}\label{lemma:ball-continuity-Phi}
	The function $\Phi:C_0([0,T],\mathbb{R}^p)\times C_0([0,T],\mathbb{R}^p)\to { \cl C_0 }^d$  is continuous on the set
	$$
		B_L=\{(f,g)\in\mathscr{H}_{(B,\hat{B})}^p:\|f\|_{H_0^{1,p}[0,T]}^2\leq L \}
	$$
for every $L>0$.
\end{lemma}
\proof Easily follows from  Lemma \ref{lemma:Phi-Phi^m-convergence} and the continuity of $\Phi_m$ (for every $m\geq 1$). \cvd

In the next lemma we will prove that $ {\mathcal J}(\cdot|(\cdot,\cdot))$ is lower semicontinuous as a function of $ ((f,g),x) \in    B_L\times { \cl C_0 }^d$.
\begin{lemma}\label{lemma:semicon-J}
	Let $((f_n,g_n),x_n)\in B_L\times { \cl C_0 }^d$ be such that
	$$((f_n,g_n),x_n)\to ((f,g),x)$$
	in $\big(C_0([0,T],\mathbb{R}^p)\times C_0([0,T],\mathbb{R}^p)\big)\times { \cl C_0 }^d$, as $n\to+\infty$.
	 Therefore,
	$$\liminf_{n\to +\infty}\mathcal{J}(x_n|(f_n,g_n))\geq  \mathcal{J}(x|(f,g)),$$
 i.e. the functional $\mathcal{J}(\,\cdot\,|(\cdot,\cdot))$  is lower semi-continuous as a function of $((f,g),x)\in B_L\times { \cl C_0 }^d$.
\end{lemma}
\proof
	By  hypothesis, $((f_n,g_n))_{n\in\mathbb{N}} \subset B_L$ is such that $(f_n,g_n)\to (f,g)$ in the space ${\cl C}_0^p\times {\cl C}_0^p$, as $n\to+\infty$. Therefore, since $B_L$ is a compact set (it is a level set of the good rate function $ I_{(B,\hat{B})}(\cdot,\cdot) $),
	we have that $(f,g)\in B_L$ and $g=\hat{f}$.
	 If $\displaystyle\liminf_{n\to +\infty}\mathcal{J}(x_n|(f_n,g_n))=\lim_{n\to +\infty}\mathcal{J}(x_n|(f_n,g_n))=+\infty$ there is nothing to prove, hence we can suppose that $(x_n)_{n\in\mathbb{N}}\subset H_0^{1,d}[0,T]$.
	 Now, combining Lemma \ref{lemma:tech2} and Remark \ref{rem:Gamma}, we have that for every $\varepsilon>0$ there exists $n_\varepsilon \in\mathbb{N}$ such that for every $n\geq n_\varepsilon$
	$$\displaylines{
		\mathcal{J}(x_n|(f_n,\hat{f}_n))=	\Gamma\Big(x_n-\Phi(f_n,\hat{f}_n)-\int_0^\cdot\mu(\hat{f}_n(t))\,dt\Big|a^{-1}(\hat{f}_n)\Big)\cr
		>(1-\varepsilon)	\Gamma\Big(x_n-\Phi(f_n,\hat{f}_n)-\int_0^\cdot\mu(\hat{f}_n(t))\,dt\Big|a^{-1}(\hat{f})\Big)\cr
		=(1-\varepsilon)\mathcal{J}\Big(x_n-\Phi(f_n,\hat{f}_n)-\int_0^\cdot\mu(\hat{f}_n(t))\,dt+\Phi(f,\hat{f})+\int_0^\cdot\mu(\hat{f}(t))\,dt\Big| (f,\hat{f})\Big)}$$
	Moreover, from Remark \ref{rem:conv-Mphi-Aphi-AphiInv}  and Lemma \ref{lemma:ball-continuity-Phi},
	$$x_n-\int_0^\cdot\mu(\hat{f}_n(t))\,dt-\Phi(f_n,\hat{f}_n)+\int_0^\cdot\mu(\hat{f}(t))\,dt+\Phi(f,\hat{f})\to x$$
	in ${ \cl C_0 }^d$, as $n\to+\infty$. 
For every $f\in H_0^{1,p}[0,T]$, the functional
	$$
		 \mathcal{J}(x|(f,\hat{f}))=\begin{cases}\displaystyle
			\Gamma\Big(x-\Phi(f,\hat{f})-\int_0^\cdot\mu(\hat{f}(t))\,dt\Big|a^{-1}(\hat{f})\Big)
			\quad&x\in H_0^{1,d}[0,T]\\
			+\infty&otherwise
		\end{cases}
	$$
	is lower semi-continuous, being
	the rate function of a LDP for the family  $(Z^{n,(f,\hat{f})})_{n\in\mathbb{N}}$, where
$Z^{n,(f,\hat{f})}(t)=\bar X^{n,\hat f}(t)+\Phi(f,\hat{f})(t)$ and $\bar X^{n,\hat f}$ is defined in (\ref{eq:barX-cond}).
 Therefore we have that
	$$\liminf_{n\to +\infty}\mathcal{J}(x_n|(f_n,\hat{f}_n))>(1-\varepsilon)\mathcal{J}(x|(f,\hat{f}))$$
	for every $\varepsilon>0$ and thus the thesis follows.
\cvd

\begin{lemma}\label{lemma:mathcalI-semicont}
	The functional  $\mathcal{I}_Z(\,\cdot\,)$  defined in (\ref{eq:I-Z-identification}) is lower semicontinuous.
\end{lemma}
\proof
	Let $L>0$ be fixed. In order to show that $\mathcal{I}_Z(\,\cdot\,)$ is lower semi-continuous, it is enough to prove that the level sets
	\begin{eqnarray*}
		M_L&=&\{x\in { \cl C_0 }^d:\mathcal{I}_Z(x)\leq L\}\\
		&=&\Bigl\{x\in H_0^{1,d}[0,T]:\inf_{f\in H_0^{1,d}[0,T]}\{I_{B,\hat{B}}(f,\hat{f})+\mathcal{J}(x|(f,\hat{f}))\}\leq L \Bigr\}
	\end{eqnarray*}
	are closed for every $L>0$. The proof is the same as Lemma 6.29 in \cite{Cel-Pac} (by using Lemma \ref{lemma:semicon-J}).
\cvd

Next Lemma is the same as Lemma 6.30 in \cite{Cel-Pac}. We give only some technical details of the proof since the multidimensional extension is not immediate.
\begin{lemma}
	If $x_m\to x$, as $m\to+\infty$, in ${ \cl C_0 }^d$, then
	$$\liminf_{m\to +\infty}I_Z^m(x_m)\geq\mathcal{I}_Z(x)$$
	where $I_Z^m(\,\cdot\,)$ and $\mathcal{I}_Z(\,\cdot\,)$ are defined in (\ref{eq:rate-function-I-Z^m}) and (\ref{eq:I-Z-identification}), respectively.
\end{lemma}
\proof If $(x_m)_{m\geq m_0}\subset { \cl C_0 }^d\backslash H_0^{1,d}[0,T]$, for some $m_0>0$
	$$\liminf_{m\to +\infty}I_Z^m(x_m)=\lim_{m\to +\infty}I_Z^m(x_m)=+\infty$$
	and there is nothing to prove, hence we can suppose that $(x_m)_{m\in\mathbb{N}}\subset H_0^{1,d}[0,T]$. Now, there are two possibilities:
	\begin{itemize}
		\item[\textit{(i)}]$\displaystyle\sup_{m\geq1}\|x_m\|_{H_0^{1,d}[0,T]}^2<+\infty$;
		\item[\textit{(ii)}]$\displaystyle\sup_{m\geq1}\|x_m\|_{H_0^{1,d}[0,T]}^2=+\infty$.
	\end{itemize}
The proof of  the case \textit{(i)}, by using Lemma \ref{lemma:mathcalI-semicont} ,is the same as in the one-dimensional case (see Lemma 6.30 in \cite{Cel-Pac}).
 Now, let us consider the case \textit{(ii)}, hence, up to a subsequence, we can suppose that
\begin{equation}\label{eq:limit-I_Z^m-infinite}
	\lim_{m\to +\infty}\|x_m\|_{H_0^{1,d}[0,T]}^2=+\infty
\end{equation}
and we have to  prove that
$\lim_{m\to +\infty}I_Z^m(x_m)=+\infty.$
For every $ u>0 $ we have,
\begin{equation}\label{eq:inf-in-out}
\begin{array}{ll}\displaystyle I_Z^m(x_m)&=\displaystyle\min\Big\{\inf_{\lVert f\rVert_{H_0^{1,p}[0,T]}^2\leq  \lVert x_m\rVert_{H_0^{1,d}[0,T]}^{2u}}\!\!\!\mathcal{H}_m((f,\hat{f}),x_m), \inf_{\lVert f\rVert_{H_0^{1,p}[0,T]}^2>  \lVert x_m\rVert_{H_0^{1,d}[0,T]}^{2u}}\!\!\!\mathcal{H}_m((f,\hat{f}),x_m))\Big\}\\
&\displaystyle\geq \min\Big\{ \inf_{\lVert f\rVert_{H_0^{1,p}[0,T]}^{2}\leq  \lVert x_m\rVert_{H_0^{1,d}[0,T]}^{2u}}\!\!\!{\mathcal J}^m(x_m|(f,\hat{f})), \inf_{\lVert f\rVert_{H_0^{1,p}[0,T]}^{2}>  \lVert x_m\rVert_{H_0^{1,d}[0,T]}^{2u}}\!\!\!I_{(B,\hat{B})}(f,\hat{f})\Big\}
\end{array}\end{equation}
where $ I_{(B,\hat{B})}(\cdot,\cdot) $ and $ \mathcal{J}^m(\cdot|(f,\hat{f})) $ are defined,  respectively, in (\ref{eq:rate-function-B-hatB}) and (\ref{eq:rate-function-J^m}). Now we consider the two infima in (\ref{eq:inf-in-out}).
For the second one we have,
\begin{equation} \label{eq:inf-out}
\inf_{\lVert f\rVert_{H_0^{1,p}[0,T]}^2>  \lVert x_m\rVert_{H_0^{1,d}[0,T]}^{2u}}I_{(B,\hat{B})}(f,\hat{f})=\inf_{\lVert f\rVert_{H_0^{1,p}[0,T]}^2> \lVert x_m\rVert_{H_0^{1,d}[0,T]}^{2u}}\frac12\lVert f\rVert_{H_0^{1,p}[0,T]}^2\geq  \frac12\,\lVert x_m\rVert_{H_0^{1,d}[0,T]}^{2u}.\end{equation}
Suppose now that $\lVert f\rVert_{H_0^{1,p}[0,T]}^{2}\leq  \lVert x_m\rVert_{H_0^{1,d}[0,T]}^{2u}$. From  Remark \ref{rem:unif-def-pos}, Remark \ref{rem:eigenvalue-estimate} and
the Cauchy-Schwarz inequality, there exists $M>0$ such that

\begin{eqnarray*}\mathcal{J}^m(x_m|(f,\hat{f}))&=&\Gamma\Big( x_m-\Phi^m(f,\hat{f}) -\int_0^\cdot\mu(\hat{f}(t))\,dt\Big|a^{-1}(\hat{f})\Big)\\
	&\geq& \frac1{2M^2\lVert f\rVert_{H_0^{1,p}[0,T]}^{2\alpha}} \int_0^T\|\dot{x}_m(t)-\dot{\Phi}^m(f,\hat{f})(t)-\mu(\hat{f}(t))\|^2\,dt\,\\
		&\geq&\frac1{2M^2\lVert f\rVert_{H_0^{1,p}[0,T]}^{2\alpha}}\int_0^T(\|\dot{x}_m(t)\|^2
		-2\|\dot{x}_m(t)\|\|\dot{\Phi}^m(f,\hat{f})(t)+\mu(\hat{f}(t))\|)\,dt.
\end{eqnarray*}

Furthermore
\begin{eqnarray*}\int_0^T\|\dot{x}_m(t)\|\|\dot{\Phi}^m(f,\hat{f})(t)+\mu(\hat{f}(t))\|\,dt &\leq & \int_0^T\|\dot{x}_m(t)\|\|\dot{\Phi}^m(f,\hat{f})(t)\|\,dt+\int_0^T\|\dot{x}_m(t)\|\|\mu(\hat{f}(t))\|\,dt
\end{eqnarray*}
and under Assumption \ref{ass:hp-mu-sigma-II}, denoting with $M$ a generic positive constant, we have that
$$\int_0^T\|\dot{\Phi}^m(f,\hat{f})(t)\|^2\,dt\leq   M \|\hat f\|^{2\alpha}_\infty \|f\|_{H_0^{1,p}[0,T]}^2\leq  M \|f\|_{H_0^{1,p}[0,T]}^{2\alpha +2},$$
and thanks to Remark \ref{rem:hat-f-bounded} and Assumption \ref{ass:hp-mu-sigma-II}  we have
$$
\int_0^T\|\mu(\hat{f}(t))\|^2\,dt\leq M \|\hat f\|^{2\alpha}_\infty\leq M \|f\|_{H_0^{1,p}[0,T]}^{2\alpha}.$$
Therefore, denoting with $M_1$ and $M_2$ two generic positive constants,
$$\int_0^T\|\dot{x}_m(t)\|\|\dot{\Phi}^m(f,\hat{f})(t)+\mu(\hat{f}(t))\|\,dt \leq
( M_1\lVert f\rVert_{H_0^{1,p}[0,T]}^{\alpha +1}+M_2 \lVert f\rVert_{H_0^{1,p}[0,T]}^{\alpha })\lVert x_m\rVert_{H_0^{1,p}[0,T]}
$$
Now we can choose $ u>0$ such that $ (\alpha+1)u<1$, and
since $\lVert f\rVert_{H_0^{1,p}[0,T]}^2\leq  \lVert x_m\rVert_{H_0^{1,d}[0,T]}^{2u} $, for large $m$ and a suitable $c>0$,
 one has
\begin{equation}\label{eq:inf-in}
\begin{array}{l}
{\mathcal J}^m(x_m|(f,\hat{f}))\geq
  \frac1{2M^2\, \lVert x_m\rVert_{H_0^{1,d}[0,T]}^{2\alpha\, u}}( \lVert x_m\rVert_{H_0^{1,d}[0,T]}^2-2M_1
  \lVert x_m\rVert_{H_0^{1,d}[0,T]}^{1+\alpha\, u}
- 2M_2\lVert x_m\rVert_{H_0^{1,d}[0,T]}^{(\alpha +1)u +1} )\\
\phantom{{\mathcal J}^m(x_m|(f,\hat{f}))}\geq c  \lVert x_m\rVert_{H_0^{1,d}[0,T]}^{2(1-\alpha u)}.
\end{array}\end{equation}
So  (\ref{eq:limit-I_Z^m-infinite}) follows from (\ref{eq:inf-in-out}), (\ref{eq:inf-out}), (\ref{eq:inf-in})  and then the proof is complete.
\cvd

Let's now check that a LDP holds for the family of log-price at final time $T$, $(Z^n(T))_{n\in\mathbb{N}}$. 
\begin{corollary}
The family of random variables $(Z^n(T))_{n\in\mathbb{N}}$ satisfies a LDP on $(\mathbb{R}^d,\mathscr{B}(\mathbb{R}^d))$ with the speed
 $\ep_n^{-2}$ and the good rate function ($z\in \R^d$),
\begin{equation}\label{eq:rate-fin}
		I_T(z)\!=\!\!\!\!\inf_{f\in H_0^{1,p}[0,T]}\biggl\{\frac12\|f\|_{H_0^{1,p}[0,T]}^2+\frac12(z-\Phi(f,\hat{f})(T)-M_{\hat{f}})^TA_{\hat{f}}^{-1}(z-\Phi(f,\hat{f})(T)-M_{\hat{f}})\biggr\},
	\end{equation}	
where
\begin{enumerate}
	\item[\textit{(i)}] the matrix $A_{\hat{f}}\in\mathbb{R}^{d\times d}$ is  such that $\displaystyle(A_{\hat{f}})_{ij}=\int_0^T a_{ij}(\hat{f}(t))\, dt$, $1\leq i,j\leq d$;
	\item[\textit{(ii)}]the vector $M_{\hat{f}}\in\mathbb{R}^d$ is such that $\displaystyle(M_{\hat{f}})_i=\int_0^T \mu_i(\hat{f}(t))\, dt$, $1\leq i\leq d$.
\end{enumerate}

\end{corollary}
\proof

 We know, from Theorem \ref{th:LDP-main}, that a LDP on ${ \cl C_0 }^d$ with the  speed $\ep_n^{-2}$ and the good rate function $I_Z(\,\cdot\,)$ holds for the family of processes $(Z^n)_{n\in\mathbb{N}}$, where $I_Z$ is the functional defined in (\ref{eq:I-Z-identification}). Let $H$ be the function defined by
$$
	H:{ \cl C_0 }^d\to\mathbb{R}^d\quad\quad
	z(\,\cdot\,)\mapsto z(T).
$$
It is not hard to prove that $H$ is a continuous function.
Therefore, by the contraction principle, 
\begin{eqnarray*}
	I_T(z)&=&\inf_{x\in H_0^{1,d}[0,T]: x(T)=z}\inf_{f\in H_0^{1,p}[0,T]}\biggl\{ \frac12 \|f\|_{H_0^{1,p}[0,T]}^2+J(x-\Phi(f,\hat{f})|\hat{f}) \biggr\}\\
	&=&\inf_{f\in H_0^{1,p}[0,T]}\inf_{x\in H_0^{1,d}[0,T]: x(T)=z}\biggl\{ \frac12 \|f\|_{H_0^{1,p}[0,T]}^2+J(x-\Phi(f,\hat{f})|\hat{f}) \biggr\}
\end{eqnarray*}
Since the term $\displaystyle\frac12\|f\|_{H_0^{1,p}[0,T]}^2$ is not dependent on $x$, we only need to calculate
\begin{equation*}
	\displaystyle\inf_{x\in H_0^{1,d}[0,T]: x(T)=z} \frac12\int_0^T(\dot{x}(t)-\mu(\hat{f}(t))-\dot\Phi(f,\hat{f})(t))^Ta^{-1}(\hat{f}(t))(\dot{x}(t)-\mu(\hat{f}(t))-\dot\Phi(f,\hat{f})(t))\,dt.
\end{equation*}
Therefore, given $z\in\mathbb{R}^d$, $f\in H_0^{1,p}[0,T]$ and $x\in H_0^{1,d}[0,T]$, set $u=\dot{x}$; we need to solve the variational calculus  problem (see for example
\cite{Gel-Sil-Fom}) with functional
\begin{equation*}
	\mathcal{F}(u)=\int_0^TF(u)\,dt=\frac12\int_0^T(u(t)-\mu(\hat{f}(t))-\dot\Phi(f,\hat{f})(t))^Ta^{-1}(\hat{f}(t))(u(t)-\mu(\hat{f}(t))-\dot \Phi(f,\hat{f})(t))\,dt
\end{equation*} and integral constraint
\begin{equation*}
	\mathcal{G}(u)=\int_0^TG(u)\,dt=\int_0^T(u(t)-\frac zT)\,dt=0
\end{equation*}
Observe that the constraint is such that
$$\int_0^Tu(t)\,dt=\int_0^T\dot{x}(t)\,dt=x(T)=z.$$
The Euler-Lagrange equation associated to the problem is
$$\frac{\partial }{\partial u}(F+\lambda^T G)=a^{-1}(\hat{f}(t))(u(t)-\mu(\hat{f}(t))-\dot \Phi(f,\hat{f})(t))+\lambda=0,$$
and then we can conclude
\cvd

\section{Short-time large deviations}\label{sect:short}
In this section we prove  a multi-dimensional short-time LDP, when $\mu=0$, following \cite{GiPaPi}.  Let us denote $\hat B^n(t)=\hat B(\delta_n t)$ and suppose that the family $(\ep_n B,\hat B^n)$ satisfies a LDP
with the speed $\ep_n^{-2}$ (see Example ...). Consider  the family of processes 
$(Z(\ep_n t))_{t\in[0,1]}$,  where $Z$ is the log-price process with $\mu=0$. We have the following result.

\begin{theorem}
In the hypotheses of Theorem \ref{th:LDP-main}
the two families   $((\ep_n \delta_n^{-1/2}Z(\delta_nt))_{t\in[0,T]})_{n\in\mathbb{N}}$ and $(Z^n)_{n\in\mathbb{N}}$  are exponentially equivalent (see Definition in DZ) and therefore satisfy the same LDP. In particular,

(i)  the family $((\ep_n \delta_n^{-1/2}Z(\delta_nt))_{t\in[0,T]})_{n\in\mathbb{N}}$ satisfies a LDP  with  speed $\ep_n^{-2}$ and good rate function given by (\ref{eq:rate-main}) with $\mu=0$;

(ii)  the family of random variables $ (\ep_n \delta_n^{-1/2}Z(\delta_n T))_{n \in \mathbb{N}}$
satisfies a LDP  with  speed $\ep_n^{-2}$ and good rate function given by (\ref{eq:rate-fin}) with $\mu=0$.
	\end{theorem}
\proof
Consider now the process 
$(Z(\ep_n t))_{t\in[0,1]}$.

Thanks to Theorem \ref{th:LDP-main} we have a LDP for the family of processes $(Z^n)_{n\in \N}$. Consider now the process 
$(Z(\delta_n t))_{t\in[0,1]}$,  In law we have,
$$\displaylines{
	 Z_i(\delta_n t)=\int_0^{\delta_n t}\bigg( -\frac12\sum_{j=1}^d\sigma_{ij}(\hat{B}(s))^2-\frac12\sum_{\ell=1}^p\tilde{\sigma}_{i\ell}(\hat{B}(s))^2 \bigg)\,ds \cr
+\int_0^{\delta_n t}\sum_{\ell=1}^p\tilde{\sigma}_{i\ell}(\hat{B}(s))\,dB_\ell(s)+\int_0^{\delta_n t}\sum_{j=1}^d\sigma_{ij}(\hat{B}(s))\,dW_j(s)\cr =
\delta_n\int_0^{ t}\bigg( -\frac12\sum_{j=1}^d\sigma_{ij}(\hat{B}^n(s))^2-\frac12\sum_{\ell=1}^p\tilde{\sigma}_{i\ell}(\hat{B}(s))^2 \bigg)\,ds \cr
+\sqrt{\delta_n}\int_0^{ t}\sum_{\ell=1}^p\tilde{\sigma}_{i\ell}(\hat{B}^n(s))\,dB_\ell(s)+\sqrt{\delta_n}\int_0^{ t}\sum_{j=1}^d\sigma_{ij}(\hat{B}^n(s))\,dW_j(s)
}$$
for every $1\leq i\leq d$ and $0\leq t \leq T$.
Now, if we repeat the proof of Theorem 2.6 in \cite{GiPaPi} we obtain that the family of processes  $((\ep_n \delta_n^{-1/2}Z(\ep_n t))_{t\in [0,T]})_{n\in \N}$
is exponentially equivalent (see Definition \ref{def:exp-equiv})   to the family $(Z^n)_{n\in \N}$ and therefore satisfy a LDP with the speed $\ep_n^{-2}$ and the rate function defined in \eqref{eq:rate-main} with $\mu=0$.
\cvd

\paragraph{Acknowledgements.}
The authors thank Paolo Pigato  for some hints and comments  about the financial aspect of the problem.
\begin{appendices}

\section{Technical results}\label{app:tech}

We collect here some  (well known) facts on continuous functions and some technical results on positive definite matrices.

\begin{remark}\label{rem:continuous}\rm
$(i)$
Suppose  $f:\R^p\to \R$ is a   continuous function and let
 $\varphi_n, \varphi\in \cl{C}^p$ be  functions
such that $ \varphi_n \overset{\cl{C}^p}{\underset{}{\longrightarrow}} \varphi, $ as $ n \to +\infty, $
then $ f\circ\varphi_n \overset{\cl{C}^p}{\underset{}{\longrightarrow}} f\circ\varphi, $ as $ n \to +\infty. $

$(ii)$
Suppose  $f:\R^p\to \R$ is a   continuous function  and let $(\varphi_n)_n\subset  C([0,T], \R^p)$ be a sequence of equi-bounded functions, i.e., there exist $M>0$ such that for every $n\in\N$,
$||\varphi_n||_\infty\leq M$, then there exist constants ${\underline f}_M,{\overline f}_M >0$ such that, for every $n\in\N$ and for every $t\in[0,T]$, $$0<{\underline f}_M\leq |f(\varphi_n(t))|\leq {\overline f}_M.$$

\end{remark}

The following properties  immediately follow  from Remark \ref{rem:continuous}, since a converging sequence of functions is an equi-bounded set in the space of continuous functions.
\begin{remark}\label{rem:conv-Mphi-Aphi-AphiInv}\rm
	Let $(\varphi_n)_{n\in\N}\subset  \cl C ^p$ and $\varphi\in  \cl C ^p$ such that $\varphi_n\to\varphi$ in $ \cl C ^p$, as $n\to+\infty$. Then,
	\begin{itemize}
		\item[\textit{(i)}]  $	\mu(\varphi_n)\to\mu(\varphi)$ and  therefore $
		\int_0^\cdot\mu(\varphi_n(s))\, ds\to \int_0^\cdot\mu(\varphi(s))\, ds
	$
	in $ \cl C^d$, as $n\to+\infty$.
	\item[\textit{(ii)}] $\sigma(\varphi_n)\to\sigma(\varphi)$ and  $a(\varphi_n)\to a(\varphi)$ in $C([0,T],\R^{d\times d})$, as $n\to+\infty$, then
	for every $1\leq i,j\leq d$, as  $n\to+\infty$,
	$\int_0^\cdot a_{ij}(\varphi_n(s))\, ds\to \int_0^\cdot a_{ij}(\varphi(s))\, ds$
	 in $C([0,T],\R)$.

	\item[\textit{(iii)}] 	Thanks to Remarks \ref{rem:unif-def-pos} and \ref{rem:continuous},
	$
		a^{-1}_{ij}(\varphi_n)\to a^{-1}_{ij}(\varphi)
	$
	in $C([0,T],\R)$, as $n\to+\infty$, for every $1\leq i,j\leq d$.
	\end{itemize}
\end{remark}

Now, let us consider the following results.
\begin{lemma}\label{lemma:tech1}
	Let $(\varphi_n)_{n\in\N}\subset  \cl C ^p$ and $\varphi\in  \cl C ^p$ such that $\varphi_n\to\varphi$ in $ \cl C ^p$, as $n\to+\infty$. Then, there exists a constant $C_{\varphi}>0$ such that
	$$v^Ta^{-1}(\varphi_n(t))v\geq C_{\varphi}>0\quad\mbox{ and }\quad v^Ta^{-1}(\varphi(t))v\geq C_{\varphi}>0$$
	for every $t\in[0,T]$, $n\in\N$ and unit vector $v\in\R^d$.

\end{lemma}
\proof Since $\varphi_n\to\varphi$ in $ \cl C ^p$, then there exists $N>0$ such that $||\varphi_n||_{\infty}\leq N$, therefore
	the proof is a consequence of Remarks \ref{rem:unif-def-pos} with
	$C_{\varphi}=\inf_{y\in [-N,N]^p}\lambda_{min}(y)>0$ and $\lambda_{min}(y)$ is the smallest eigenvalue of $a^{-1}(y)$.
\cvd
\begin{lemma} \label{lemma:tech2}
	Let $(\varphi_n)_{n\in\N}\subset  \cl C ^p$ and $\varphi\in  \cl C ^p$ such that $\varphi_n\to\varphi$ in $ \cl C ^p$, as $n\to+\infty$. Then, for every $\varepsilon>0$ there exists $n_{\varepsilon}\in\N$ such that the matrix $a^{-1}(\varphi_n(t))-(1-\varepsilon)a^{-1}(\varphi(t))$ is (strictly) positive definite, for every $n\geq  n_\varepsilon$ and $0\leq t\leq T$.
\end{lemma}
\proof
	Fix $\varepsilon>0$ and a unit vector $v\in\R^d$. Then from Lemma \ref{lemma:tech1}, for every $t\in[0,T]$, we have
	\begin{eqnarray*}
		v^T(a^{-1}(\varphi_n(t))-(1-\varepsilon)a^{-1}(\varphi(t)))v&=&	v^T(a^{-1}(\varphi_n(t))-a^{-1}(\varphi(t)))v+\varepsilon v^Ta^{-1}(\varphi(t))v\\
		&\geq&v^T(a^{-1}(\varphi_n(t))-a^{-1}(\varphi(t)))v+\varepsilon C_{\varphi}
	\end{eqnarray*}
Now, from \textit{(iii)} in Remark \ref{rem:conv-Mphi-Aphi-AphiInv},  there exists $n_{\varepsilon}\in\N$ such that for every $n\geq n_{\varepsilon}$ and  $t\in[0,T]$ one has
$$|v^T(a^{-1}(\varphi_n(t))-a^{-1}(\varphi(t)))v|\leq \frac{\varepsilon C_{\varphi}}{2}$$
Therefore, for every $t\in[0,T]$ and $n\geq n_{\varepsilon}$
$$v^T(a^{-1}(\varphi_n(t))-(1-\varepsilon)a^{-1}(\varphi(t)))v\geq\frac{\varepsilon C_{\varphi}}{2}>0.\quad \quad \mbox{\cvd}$$

\begin{lemma}\label{lemma:tech3}
	Let $(\varphi_n)_{n\in\N}\subset  \cl C ^p$ and $\varphi\in  \cl C ^p$ such that $\varphi_n\to\varphi$ in $ \cl C ^p$, as $n\to+\infty$. Then, there exists $M>1$ such that the matrix $Ma^{-1}(\varphi_n(t))-a^{-1}(\varphi(t))$ is strictly positive definite, for every $n\in\N$  and $0\leq t\leq T$.
\end{lemma}
\proof
	Let $v\in\R^d$ be a unit vector and $M>1$. From Lemma \ref{lemma:tech1} for every $t\in[0,T]$,
	\begin{eqnarray*}
		v^T(Ma^{-1}(\varphi_n(t))-a^{-1}(\varphi(t)))v&=&(M-1)v^Ta^{-1}(\varphi_n(t))v+v^T(a^{-1}(\varphi_n(t))-a^{-1}(\varphi(t)))v\\
		&\geq&(M-1)C_{\varphi}+v^T(a^{-1}(\varphi_n(t))-a^{-1}(\varphi(t)))v.
	\end{eqnarray*}
 From \textit{(iii)} in Remark \ref{rem:conv-Mphi-Aphi-AphiInv}, there exists a constant $N>0$ such that for every $n\in\N$
$$|v^T(a^{-1}(\varphi_n(t))-a^{-1}(\varphi(t)))v|\leq N$$
 Therefore, for every $t\in[0,T]$ and $n\in\N$,
$$v^T(Ma^{-1}(\varphi_n(t))-a^{-1}(\varphi(t)))v\geq (M-1) C_{\varphi}-N$$
thus, it is enough to choose $M>1$  such that $(M-1)C_{\varphi}-N>0.$
\cvd

\begin{remark}\label{rem:hat-f-bounded}\rm
For $L>0$, denote by $D_{L}$ the closed ball in the Cameron Martin space, i.e.
$$
D_{L}=\{f\in H_0^{1,p}[0,T]: \lVert f\rVert^2_{H_0^{1,p}[0,T]}\leq L\}.$$
Then for $f\in D_{L}$, from the Cauchy-Schwarz inequality,
$$||\hat{f}(t)||^2=\sum_{\ell=1}^p\hat{f}_\ell(t)^2=\sum_{\ell=1}^p\bigg( \int_0^TK_\ell(t,s)\dot{f}_\ell(s)\,ds \bigg)^2\leq \sum_{\ell=1}^p\int_0^TK_\ell(t,s)^2\,ds\int_0^T\dot{f}_\ell(t)^2\,ds. $$
	Therefore (thanks to condition $(b)$  in Definition \ref{def:Volterra-process}) there exists a constant $ M>0 $
such that  $$||\hat{f}(t)||^2\leq M ||{f}||_{H_0^{1,p}[0,T]}^2$$ and therefore
$$\sup_{f\in D_{L}}\sup_{t \in [0,T]}||\hat{f}(t)||\leq M L.$$
\end{remark}
\end{appendices}

\end{document}